\documentclass[pdflatex,sn-mathphys-num]{sn-jnl}
\usepackage{amsmath,amssymb, amsthm}
\usepackage{lipsum}  
\usepackage{dsfont}
\usepackage{subcaption}	
\usepackage{enumerate}
\usepackage{hyperref}
\usepackage{multicol}
\usepackage[dvipsnames]{xcolor}
\usepackage{hyperref}
\usepackage[numbers]{natbib}
\usepackage{MnSymbol} 

\setcitestyle{square,  numbers}

\usepackage{float}
\usepackage{graphicx}
\graphicspath{ {./} }

\newcommand{\tgeq}[1]{(#1)_{t \geq 0}}
\newcommand{\tinT}[1]{(#1)_{t \in [0,T]}}
\newcommand{\mcal}[1]{\mathcal{#1}}
\newcommand{\mbb}[1]{\mathbb{#1}}

\providecommand{\keywords}[1]{\textbf{Keywords and phrases: } #1}

\newcommand{\R}{\mbb{R}}
\newcommand{\N}{\mbb{N}}

\newcommand{\tr}{\textnormal{tr}}

\newcommand{\dif}{\textnormal{d}}

\newcommand{\1}{\mathds{1}}

\newcommand{\Borel}{\mcal{B}}

\newcommand{\gen}{\mcal{L}}

\newcommand{\Prob}{\mbb{P}}
\newcommand{\Qrob}{\mbb{Q}}

\newcommand{\E}{\mbb{E}}

\newcommand{\cond}{\, | \,}

\newcommand{\filtt}{\mathcal{F}_t}
\newcommand{\Yfiltt}{\mathcal{Y}_t}

\newcommand{\markspace}{\mathcal{K}}

\makeatletter
\newcommand{\customlabel}[2]{%
  \def\@currentlabel{#2}%
  \label{#1}%
  #2%
}
\makeatother

\usepackage{amsthm}

\theoremstyle{plain}
\newtheorem{assumption}{Assumption}

\newtheorem{theorem}{Theorem}[section] 
\newtheorem{lemma}[theorem]{Lemma}
\newtheorem{proposition}[theorem]{Proposition}
\newtheorem{corollary}[theorem]{Corollary}

\newtheorem{definition}[theorem]{Definition} 
\newtheorem{example}[theorem]{Example}
\newtheorem{remark}[theorem]{Remark}

\begin{document}

\title[Nonlinear Filtering and Spatial Asymptotic Consistency for SPDEs]{Nonlinear Filtering and Spatial Asymptotic Consistency for SPDEs Observed via Spatio-Temporal Point Processes}

\author*[1]{\fnm{Jan} \sur{Szalankiewicz}}\email{szalankiewicz@math.tu-berlin.de}
\author[2]{\fnm{Cristina} \sur{Martinez-Torres}}\email{martineztorres@uni-potsdam.de}
\author[1]{\fnm{Wilhelm} \sur{Stannat}}\email{stannat@math.tu-berlin.de}

\affil*[1]{\orgdiv{Institute of Mathematics}, \orgname{TU Berlin}, \orgaddress{\postcode{10623}, \city{Berlin}, \country{Germany}}}
\affil[2]{\orgdiv{Institute of Physics and Astronomy}, \orgname{University of Potsdam}, \orgaddress{\postcode{14476}, \city{Potsdam}, \country{Germany}}}

\abstract{In this paper, we develop the mathematical framework for filtering problems arising from biophysical applications where data is collected from confocal laser scanning microscopy recordings of the space-time evolution of intracellular wave dynamics of biophysical quantities. In these applications, signals are described by stochastic partial differential equations (SPDEs) and observations can be modelled as functionals of marked point processes whose intensities depend on the underlying signal. We derive both the unnormalized and normalized filtering equations for these systems, demonstrate the asymptotic consistency and approximations of finite dimensional observation schemes respectively partial observations. Our theoretical results are validated through extensive simulations using synthetic and real data. These findings contribute to a deeper understanding of filtering with point process observations and provide a robust framework for future research in this area.}


\maketitle

\section{Introduction}

Reaction-diffusion systems are fundamental models in biophysics, representing spatially extended systems where dynamics at each location involve nonlinear reaction kinetics, coupled by diffusive transport of reacting species \cite{source_spde_life_sciences_1, source_spde_life_sciences_2}. The 
motivating example for this paper is the spatially extended stochastic FitzHugh-Nagumo-type model of actin wave formation in the social amoeba \textit{Dictyostelium discoideum} \cite{source_dicty_model}, modeled by a stochastic partial differential equation (SPDE) of the following type:
\begin{equation}\label{introduction_spde}
\quad
    \begin{cases}
    \begin{aligned}
            \dif X(t,x) &= (AX(t,x) +F(X(t,x)) ) \dif t + B(X(t,x)) \dif W(t,x),  \\
            X(0,x) &= \xi(x),
    \end{aligned}
    \end{cases}
\end{equation}
$ t \in (0,T]$,  on a suitable domain $\mcal{D}\subset \R^d$, where $A$ denotes diffusion, and $F$ the reaction-kinetics; see \cite{source_dicty_spde_1}. We will give precise conditions on the above terms  in Section \ref{subsection_spde_solutions}. \newline

In practice, information on actin wave dynamics is obtained from confocal laser scanning microscopy (CLSM) recordings given as a time 
series of digital grey-scale images. To infer accurate statistical information contained in the data,  
based on the SPDE model \eqref{introduction_spde}, requires careful selection of the model parameters 
guided by experimental data obtained from CLSM recordings of giant 
\textit{D. discoideum} cells. In addition, even if detailed simulations based on 
\eqref{introduction_spde} may align well with experimental data, questions about the robustness and 
plausibility of model parameters remain \cite{source_dicty_spde_1}. \newline

In order to gather data using CLSM in the experiments, cells are tagged with fluorescent biomarkers, allowing researchers to count photon emissions correlated with the actin 
concentration rather than measuring actin concentration directly. Each pixel in the CLSM images 
corresponds to a specific region of the cell, with pixel values representing the number of emitted photons. 
Consequently, CLSM recordings provide data as sequences of digital images, where the photon counts are 
approximately Poisson distributed with intensity related to the fluorescent material concentration. This 
introduces an additional layer of stochasticity known as \textit{observation noise}. 
 \newline

We use marked point processes (MPPs) as a mathematical model of such type of observations. 
MPPs represent a well-established class of point processes, capable of modeling random events in random 
positions — in this case, the time and location of photon emissions. This approach allows us to infer 
information on the underlying signal, the actin concentration modelled in terms of the SPDE 
\eqref{introduction_spde}, given MPP observations using stochastic filtering, a comprehensive Bayesian 
framework for sequential estimation in a model-based setting.  \newline

More specifically, let $\markspace$ be the \textit{mark space} modeling the area of point positions, the 
evolution of the photon emissions in a given subset $\Gamma \subseteq \markspace$ over time can be 
written in integral form as the dynamics of a stochastic jump process $Y$ as follows,
\begin{equation} \label{introduction_MPP}
\quad
    \begin{cases}
    \begin{aligned}
             \dif Y(\Gamma,t) &= \int_\Gamma \lambda(t ,x \cond X(t)) \dif x \, \dif t + \dif N(\Gamma,t) , \quad t \in (0,T], \\
             Y(\Gamma,0) &= 0,
    \end{aligned}
    \end{cases}
\end{equation} 
where $ \tgeq{N(\Gamma,t)}$ is the jump martingale corresponding to $Y$ restricted to $\Gamma$. \newline

In this paper, we develop the statistical filtering theory for the stochastic signal $X$ described by the SPDE in \eqref{introduction_spde} with observation schemes arising from \eqref{introduction_MPP}. Our work includes the derivation of the Kallianpur-Striebel formula, as well as the Zakai and Kushner-Stratonovich equations for the posterior distribution of $X$. Although filtering problems are often formulated with Gaussian observations \cite{filtering_fundamentals}, the study of filtering with point process observations has gained significant attention across various disciplines, including statistics and engineering \cite{source_filtering_point_processes_3, source_filtering_point_processes_2, source_filtering_point_processes_1}.

The foundational work by Snyder \cite{Snyder_filtering} was the first to rigorously address point process 
observations in stochastic filtering, a framework later extended to MPPs by Brémaud 
\cite{Bremaud_phd, Bremaud}. Filtering for SPDEs with Gaussian observations was initially explored by Pardoux \cite{Pardoux_SPDE_filter}, and further developed by Ahmed, Fuhrmann, and Zabczyk \cite{filtering_Ahmed_Fuhrmann_Zabczyk}. Florchinger made contributions by analyzing SPDE signals with one-dimensional temporal point process observations \cite{Florchinger_filtering_spde}, though this line of inquiry was not extensively pursued. More recently, Sun, Zeng, and Zhang investigated filtering with MPPs in the context of abstract Hilbert-space valued Markov processes \cite{Wei_Sun_filtering_MPP}, albeit without deriving the Kushner-Stratonovich equation and without giving an explicit functional analytical framework for the signal process.

To the best of our knowledge, the filtering framework of SPDEs with multivariate point process 
observations or more general MPP observations, has not been previously addressed in the 
literature. \newline

Furthermore, we explore the relationship between observations represented as marked point processes and 
their lower-resolution multivariate point process approximations, which contain reduced spatial 
information. We prove weak convergence of the multivariate point processes observations to the 
underlying MPP counterparts and establish convergence in total variation for both, the unnormalized and 
normalized posterior distributions in the high-resolution limit. Additionally, we address the case of 
partial observations. To the best of our knowledge, such work has not yet been conducted within the 
context of filtering, providing error bounds for estimates based on low-resolution point process 
observations. Finally, we report on extensive numerical experiments, providing further insights into our 
theoretical findings. \newline

The structure of this paper is as follows. In Section \ref{sec_2} we provide a concise overview of key 
concepts of SPDEs in the variational setting and MPPs, followed by the precise mathematical modeling 
of the stochastic filtering problems including both infinite- and finite-dimensional spatio-temporal 
point process observation schemes. 

Section \ref{section_filtering_equations} is devoted to deriving the filtering equations. Specifically, 
we present the Kallianpur-Striebel formulas in Lemma \ref{lemma_kallianpur_striebel} and Lemma 
\ref{lemma_kallianpur_striebel_multivariate}, the Zakai equations for the time-evolution of the 
unnormalized conditional distributions in Theorem \ref{Zakai_thm} and Theorem 
\ref{Zakai_thm_multivariate} and the Kushner-Stratonovich equations in Corollary 
\ref{kushner_strat_thm_mpp} and Corollary \ref{kushner_strat_thm_multivariate}.

In Section \ref{section_observation_schemes}, we study the convergence of the multivariate point 
processes observations to the underlying MPP counterparts in the high-resolution limit,  
analyze the convergence of both, the unnormalized and normalized posterior distributions and establish 
approximation errors. Additionally, we introduce a specific modeling of partial observations designed to 
replicate the setting of CLSM data and derive corresponding error bounds.

The final Section \ref{sec_simulations} presents numerical simulation results.

\section{Mathematical setting of the filtering model}\label{sec_2}

The filtering theory for SPDE signals with Gaussian observations has been extensively studied in the literature; see \cite{Pardoux_SPDE_filter, filtering_Ahmed_Fuhrmann_Zabczyk}. The only known work analyzing SPDE signals with point process observations is the conference paper \cite{Florchinger_filtering_spde}, which considers a one-dimensional Poisson process with intensity dependent on the SPDE state. The recent paper \cite{Wei_Sun_filtering_MPP} introduces multivariate point process (MPP) observations but deals with a very abstract, Hilbert space-valued Markov process.

Our objective to explicitly model the CLSM observations of actin wave dynamics implies leads to a new filtering problem for an SPDEs observed with MPPs. For one, this approach introduces a novel method for modeling spatio-temporal shot noise via generalized Cox processes steered by an SPDE. Furthermore, new questions about limits of statistical estimators arise, which we partly answer in Section \ref{section_observation_schemes}.
\subsection{The signal process} \label{subsec_signal_process}
We will model the signal process as an SPDE within the \textit{variational} framework as introduced in \cite{Liu_Roeckner, Pardoux_SPDE_filter}, employing their terminology.  Although our analysis primarily focuses on the variational solution concept, it can be adapted to accommodate other concepts, such as mild solutions. This adaptation is a technical matter that necessitates changes to the functional analytical framework, resulting in different conditions for the SPDE coefficients and a different Itô formula than the one we employ; see for example \cite[Thm. 4.17]{da_prato_zabczyk}.

\subsubsection{Variational solutions to SPDE} \label{subsection_spde_solutions}
Let $\mcal{H}$ be a Hilbert space with inner product $(\cdot, \cdot )_\mcal{H} $ and $\mcal{V}$ a reflexive Banach space, both on $\mcal{D}\subset \R^d $, and let $\mcal{V}^\ast$ denote the dual space of $\mcal{V}$. By $_{\mcal{V}^\ast}\langle \cdot , \cdot \rangle_\mcal{V}$ we denote the dual pairing between $\mcal{V}$ and  $\mcal{V}^\ast$. We impose that $(\mcal{V},\mcal{H}, \mcal{V}^\ast )$  forms a \textit{Gelfand triple} which implies that $\mcal{V} \subset \mcal{H} \approx \mcal{H}^\ast \subset \mcal{V}^\ast $ continuously and densely and that 
\begin{equation*}
    _{\mcal{V}^\ast}\langle h , v  \rangle_\mcal{V} = (h, v )_\mcal{H}, \quad \text{ for all } h \in \mcal{H}, \, v \in \mcal{V},
\end{equation*}
see e.g. \cite[pp. 69]{Liu_Roeckner}. \\

Let $T \geq 0$ and $(\Omega, \mcal{F}, \tgeq{\filtt}, \Prob)$ be a complete probability space with filtration $\tgeq{\filtt}$ satisfying the usual conditions. For some given separable real Hilbert space $\mcal{U}$ we consider $(W(t))_{t \geq 0}$ to be a $\mcal{U}$-valued $\tgeq{\filtt}$-adapted $Q$-Wiener process.  We assume that $Q$  is a self-adjoint, positive semidefinite linear operator on $\mcal{U}$, with finite trace $\textnormal{tr}_\mcal{U} Q < +\infty$.

We consider stochastic partial differential equations (SPDE) on $\mcal{H}$ of the following type
\begin{equation*} \customlabel{introduction_spde_basic_formulation}{\textnormal{(S)}}
\quad
    \begin{cases}
    \begin{aligned}
            \dif X(t) &= A(X(t)) \dif t + B(X(t)) \dif W(t), \quad t \in (0,T], \\
            X(0) &= \xi\in \mcal{H}
    \end{aligned}
    \end{cases}
\end{equation*}
with $B \in L_2(\mcal{U}, \mcal{H})$, where $L_2(\mcal{U}, \mcal{H})$ denotes the space of Hilbert Schmidt operators from $\mcal{U}$ to $\mcal{H}$, and $A: \mcal{V} \rightarrow \mcal{V}^\ast$. Such a general form of an SPDE covers cases such as stochastic heat and reaction-diffusion equations, see \cite{Liu_Roeckner}. In order to being able to work with an analytically weak solution to \ref{introduction_spde_basic_formulation} we make the standard assumptions:

\begin{assumption} \label{ass_signal} We assume that the following conditions hold on the coefficients $\xi, \, A, B$ in \ref{introduction_spde_basic_formulation}. \

    \begin{enumerate}
    \item[\customlabel{assumption_initial_cond}{\textnormal{(A0)}} ]\textnormal{Initial condition}: Let $\xi \in L^2(\Omega, \mcal{F}_0, \Prob; \mcal{H})$.
    \item[\customlabel{assumption_hemicts}{\textnormal{(A1)}} ]\textnormal{Hemicontinuity}: For $u, v, w \in \mcal{V}$, $t \in [0,T]$ the map
    \begin{equation*}
        \delta \mapsto \, _{\mcal{V}^\ast}\langle A(u+\delta v) , w  \rangle_\mcal{V}
    \end{equation*}
    is continuous.
    \item[\customlabel{assumption_weak_mon}{\textnormal{(A2)}}] \textnormal{Weak monotonicity}: There exists a constant $\mcal{C}_1 \in \R$ s.t. for $u,v \in \mcal{V}$
    \begin{align*}
        2 \, _{\mcal{V}^\ast}\langle A(u) - A(v) , u-v  \rangle_\mcal{V} &+ \|( B(u) - B(v) ) \, \sqrt{Q}  \|^2_{L_2(\mcal{U}, \mcal{H})} \leq \mcal{C}_1  \| u-v \|^2_\mcal{H} 
    \end{align*}
    on $[0,T]$.
    \item[\customlabel{assumption_coercivity}{\textnormal{(A3)}}] \textnormal{Coercivity}: There exist constants $\mcal{C}_2 \in \R$, $\mcal{C}_3, \mcal{C}_4 \in (1, \infty)$, $\tilde{p} \in (1,\infty)$,  such that for all $v \in \mcal{V}$
    \begin{align*}
        2 \, _{\mcal{V}^\ast}\langle A(v) , v  \rangle_\mcal{V} &+ \|  B(v)   \, \sqrt{Q} \|^2_{L_2(\mcal{U}, \mcal{H})} \leq \mcal{C}_2 \| v \|^2_\mcal{H} - \mcal{C}_3 \| v \|^{\tilde{p}}_\mcal{V} + \mcal{C}_4.
    \end{align*}
    \item[\customlabel{assumption_boundedness}{\textnormal{(A4)}}] \textnormal{Boundedness}: There exists a constant $\mcal{C}_5 >0$ s.t. for all $v \in \mcal{V}$
    \begin{align*}
        \| A(v) \|_{\mcal{V}^\ast} \leq \mcal{C}_5 (1+\|v\|_{\mcal{V}}).
    \end{align*}
\end{enumerate}
\end{assumption}

Under  Assumption \ref{ass_signal} it is known that equation \ref{introduction_spde_basic_formulation} admits an \textit{analytically weak} or \textit{variatonal solution} to the SPDE \ref{introduction_spde_basic_formulation}, see for example \cite[Thm. 4.2.4]{Liu_Roeckner}. In particular, this means that there  exists a  unique $\mcal{H}$-valued, $(\filtt)$-adapted process $X=\tinT{X(t)}$, where 
\[
X \in L^2([0,T]\times \Omega, \dif t \otimes \Prob; \mcal{H}) \cap L^{\tilde{p}}([0,T]\times \Omega, \dif t \otimes \Prob; \mcal{V}) 
\]
with $\tilde{p}$ from \ref{assumption_coercivity}, such that for any $v \in \mcal{V}$  we have the $\Prob$-a.s. equality
\begin{equation}
    (X(t), v)_\mcal{H} =  (X(0), v)_\mcal{H} + \int_0^t \hspace{-4pt}  \,   _{\mcal{V}^\ast}\langle A(X(s)) , v  \rangle_\mcal{V} \dif s + \int_0^t (v, B(X(s))\dif W(s))_\mcal{H},
\end{equation}
for any $t \in [0,T]$. Additionally, one can show that the solution is an $\mcal{H}$-Markov process \cite[Proposition 4.3.5]{Liu_Roeckner}. Such a variational solution to \ref{introduction_spde_basic_formulation} represents the signal in our filtering problem.
\subsubsection{It\^{o} functions and the infinitesimal generator}
For deriving the filtering equations in Section \ref{section_filtering_equations}, it will be of great use to have a version of It\^{o}'s lemma for variational solutions. The suitable function class is given as follows.

\begin{definition}\cite[p. 136]{Pardoux_SPDE_filter} \label{def_ito_function}
  We call a function $\psi: \mcal{H} \rightarrow \R$ an \textit{It\^{o} function}, if it fulfills the following conditions, where all derivatives have to be understood w.r.t. $\mcal{H}$. \begin{enumerate}[(i)]
      \item $\psi$ is twice Fréchet-differentiable with derivatives $\text{D}^1 \psi$ and $\text{D}^2 \psi$.
      \item $\psi$,  $\text{D}^1 \psi$ and $\text{D}^2 \psi$ are locally bounded.
      \item For any trace-class operator $\Theta:\mcal{H}\rightarrow \mcal{H}$, the functional $u \rightarrow \textnormal{tr}\big( \Theta \text{D}^2 \psi(u) \big)$ is continuous on $\mcal{H}$.
      \item For $v \in \mcal{V}$ both $\text{D}^1 \psi(v) \in \mcal{V}$ and the map $\text{D}^1 \psi(v)|_\mcal{V}:\mcal{V} \rightarrow \mcal{V}$ is continuous when the domain is equipped with the strong and the image is equipped with the weak topology.
      \item There is a constant $\mcal{C}_\mcal{V} > 0$ such that $\| \text{D}^1 \psi(v) \|_\mcal{V} \leq \mcal{C}_\mcal{V}(1+\|v\|_\mcal{V})$ for all $v \in \mcal{V}$.
  \end{enumerate}

  Moreover, if $\psi$, $\text{D}^1 \psi$ and $\text{D}^2 \psi$ are globally bounded, we call $\psi$ a \textit{globally bounded It\^{o} function}. 
\flushright $\diamond$
\end{definition}

\noindent Under Assumptions \ref{ass_signal}, the \textit{infinitesimal generator} $\gen$ of the signal $X$ is given by
\begin{equation} \label{signal_generator}
    \gen \psi = \,  _{\mcal{V}^\ast}\langle A(\cdot ) , \text{D}^1 \psi   \rangle_\mcal{V} +  \frac{1}{2} \tr \{ \text{D}^2 \psi \; B(\cdot) QB(\cdot)^\ast \},
\end{equation}
for any It\^{o} function $\psi$.

\subsection{The observation process} \label{subsec_observation}
In the biophysical application we can only measure the actin concentration indirectly in the form of photon emissions of certain fluorescent biomarkers attached to actin. These measurements are given as sequences of digital gray-scale images in given times $t_1, \dots, t_n$. In particular, the pixel value of an image in time $t_{i}$ corresponds to a (transformed) photon emission count in the corresponding area under the microscope, recorded in the time interval $(t_{i-1}, t_{i}]$. In practice, our analysis has shown that in the given experiments these photon counts have a Poisson statistic. Hence, we can justify to model the photon count of an individual pixel as a Poisson distributed random variable, where the intensity is given as a function of the concentration of fluorescent molecules available at the time of recording in the corresponding area. \newline

Now, for a sequence of images, an intuitive approach to modeling such an observation scheme is to assign a point process in time to each pixel, resulting in a multivariate point process as described in \ref{introduction_MPP_discretized_observation_multivariate}, where $M$ is the number of pixels. This is referred to as the \textit{finite dimensional} model because it only involves a finite number of sets, or pixels.

A more general approach is to move away from the analogy of digital images with a fixed number of pixels and instead look at (theoretical) recordings of the exact space-time locations of each single photon count. An analytically manageable way to formalize such an observation is by employing the notion of marked point processes, which can be either seen as random space-time point clouds or as random space-time counting measures. This leads to the scheme \ref{introduction_MPP_observation}, termed the \textit{infinite dimensional} observation.

We choose to first construct the more general version \ref{introduction_MPP_observation}, as this observation includes the exact times and locations of photon emissions. From this, we derive \ref{introduction_MPP_discretized_observation_multivariate}, which records only the pixel area of photon emissions, not their exact positions. This distinction will become clearer once all technical details are elaborated. \newline

In the first half of this section, we provide a brief overview of point process theory, as several of the tools discussed are crucial for the analyses in Sections \ref{section_filtering_equations} and \ref{section_observation_schemes}. The second half introduces the two observation schemes we intend to investigate and outlines the filtering problem.

For a comprehensive introduction to point processes, we refer the reader to \cite{Bremaud, DVJ_vol2, Lipster_shiryaev}, which serve as our primary references regarding MPPs.

%
%
%
%
%
%
%

%
%
%
%
%
%
%
%
%
%
%
%
\subsubsection{Observation schemes} \label{subsection_observation_structures}

For details on point process theory we refer the reader to the Appendix. Let $(\Omega, \mcal{F},( \filtt), \Prob)$ be the filtered probability space  and $X$ be the signal from Section \ref{subsection_spde_solutions}. We impose the following assumptions.
\begin{assumption}  \label{ass_observation}
        Let $\markspace \subset \R^{d_O}$ be  compact.  The observation process $Y$ is given as a generalized Cox process on $[0,T]\times \markspace$ directed by $X$, with boundedly finite first moment. Moreover, the conditional $(\Prob, \filtt)$-intensity $\lambda$ of $Y$ is a strictly positive, bounded, measurable mapping  $\lambda:[0,T]\times \markspace  \times \mcal{H} \longrightarrow \R_+$  such that there exist constants $\underline{C} , \, \overline{C}$ with
    \begin{equation} \label{condition_boundedness_intensity}
       0< \underline{C} \leq \int_\markspace \lambda(t,x \cond u )\, \mu_\markspace(\dif x)  \leq \overline{C} <\infty, \quad \Prob\text{-a.s.}, \; u \in \mcal{H}, \; t \in [0,T].
    \end{equation}
     \flushright $\diamond$
\end{assumption}
 As discussed in the previous section, with fixed $T>0$,
 an MPP on $[0,T]\times \markspace$ is not only $\Prob$-a.s. boundedly finite but even $\Prob$-a.s. totally finite. Therefore, assuming the boundedness of the stochastic intensity is not overly restrictive in this context.

\begin{remark}
 Using the notion of local characteristics introduced in the last section, condition \eqref{condition_boundedness_intensity} is equivalent to saying that the $(\Prob, \filtt)$-local characteristics $(\lambda^g(t,X(t), \Phi(\dif x\cond t, X(t))$ of $Y$ are uniformly bounded, $\mcal{H}$-measurable mappings such that
\begin{align}
           0< &\underline{C} \leq \lambda^g(t\cond X(t))   \leq  \overline{C} <\infty, \\
           &\int_\markspace  \Phi(\dif x\cond t, X(t)) =1 
\end{align}
 $\Prob$-a.s., for any $t \in [0,T]$.   
     \flushright $\diamond$
\end{remark}

\paragraph{Infinite-dimensional observations}  Given $Y$ as in Assumption  \ref{ass_observation}, the  observation \ref{introduction_MPP_observation} is a realization of the MPP $Y$ on $[0,T] \times \markspace$ given a signal path of $X$, meaning that for any Borel set $\Gamma \in \Borel(\markspace)$, by using the form of the semimartingale decomposition in \eqref{semimartingale_decomposition_MPP_compensator}, we have a path of the jump process
\begin{equation*} \customlabel{introduction_MPP_observation}{\text{(O)}}
\quad
    \begin{cases}
    \begin{aligned}
            \dif Y_{\Gamma}(t) &= [ \lambda^g(t \cond X(t))  \Phi(\Gamma \cond t, X(t))] \, \dif t + \dif N_{\Gamma}(t) , \quad t \in (0,T], \\
             Y_{\Gamma}(0) &= 0.
    \end{aligned}
    \end{cases}
\end{equation*}

\paragraph{Finite-dimensional observations}
In \ref{introduction_MPP_observation}, for any $t \in [0,T]$, given $X$ the observation $Y_\cdot(t)$ is a measure on $(\markspace, \Borel(\markspace))$. In practice we often have a finite-dimensional observation vector, think of pixels in an image from fluorescence microscopy, which dictates a specific partition on the mark space $\markspace$, thereby limiting the available spatial information and hence the choice of test sets. A mathematical formalization of such a spatial discretization can be done as follows: For any $M \in \N$ we denote by 
\[
\markspace^M := \{ K^M_1, \dots, K^M_M \} 
\] 
a partition consisting of nonempty Borel sets of the markspace $\markspace$. Such a collection of sets $\markspace^M$ can always be found for any $M \in \N$ due to the separability assumption on $\markspace$. \newline

Given any partition $\markspace^M$ and a realization of the signal $X$, we define
\[\lambda^M_{i}(t\cond X(t)):= \lambda^g(t \cond X(t))  \Phi(K^M_i \cond t, X(t)), \quad i = 1,\dots , M,\]
for any $t \in [0,T]$. We now introduce a multivariate $M$-dimensional point process $\tinT{Y^M(t)}$ on $[0,T]$, with $Y^M(t) := (Y^M_1 (t), \dots, Y^M_M(t))$, $t \in [0,T]$, where each of the  $Y^M_i$ has $(\Prob, \filtt)$-intensity $\lambda^M_{i}(t\cond X(t))$. Exactly as in \eqref{semimartingale_decomposition_MPP_compensator}, any of the processes $Y^M_i$ can be written as a semimartingale with associated jump martingale part 
\[
     \dif N^M_i(t)  := \dif Y_i(t) - \lambda^M_i(t)  \dif t  .
\]

The finite-dimensional observation is then given as the system
    \begin{equation*} \customlabel{introduction_MPP_discretized_observation_multivariate}{\textnormal{(O$^M$)}}
    \quad
        \begin{cases}
        \begin{aligned}
                \dif Y^M_{i}(t) &= \lambda^M_i(t\cond X(t)) \, \dif t + \dif N^M_i(t) , \quad t \in (0,T], \\
                 Y^M_{i}(0) &= 0, 
        \end{aligned}
        \end{cases}
    \end{equation*}
for $i = 1,\dots , M$.

\begin{remark}
   Although we could also introduce a general multivariate point process in the form of \ref{introduction_MPP_discretized_observation_multivariate}, we choose to explicitly construct the finite dimensional observation from the MPP as this approach allows us to utilize the more general  methods in both settings from the outset. Moreover, we do not need to introduce additional assumptions on the multivariate point process $Y^M$ as they carry over from the properties of $Y$. We will furthermore have the advantage of being able to embed the multivariate point process $Y^M$ on $[0,T]$ into the space of counting measures on $[0,T]\times \markspace$ in Section \ref{section_observation_schemes}. This way we characterize $Y$ as a weak limit of multivariate point processes and show how the filtering equations for \ref{introduction_MPP_observation} can be seen as the limit case of the ones corresponding to \ref{introduction_MPP_discretized_observation_multivariate}.
        \flushright $\diamond$
\end{remark}

We end this section with a simple practical example of our observation schemes.

\begin{example}[Reaction-Diffusion SPDE with Marked Cox process observations] \label{example_Cox_process} 
For some given  bounded compact domain $\mcal{D} \subset \R^d $ and a  globally Lipschitz continuous and bounded function $F$, we define the  $A(u) := \Delta u + F(u)$ with Dirichlet boundary conditions, such that \ref{introduction_spde_basic_formulation} becomes
\begin{equation} \label{example_spde_porous}
    \dif X(t) = (\Delta X(t)  + F(X(t)))\dif t + B\dif W(t),
\end{equation}
which represents a typical reaction-diffusion SPDE. We choose $\mcal{V}:= W^{1,2}_0(\mcal{D})$, $\mcal{H} := L^2(\mcal{D})$, so $\mcal{V}^\ast:= (W^{1,2}_0(\mcal{D}))^\ast$; see \cite[Ch. 4.1]{Liu_Roeckner} for a detailed discussion. \newline

Now, we explicitly construct a simple example for a marked Cox process observation of $X$. To this end, let $\markspace = \mcal{D}$ and let $0< c_1 < c_2$. We define 
\begin{equation*}
    \lambda^g(u) := \max\big\{\| u \|_\mcal{H} + c_1, c_2 \big\}, \quad u \in \mcal{H}.
\end{equation*}
For some given mollifier $\varphi_\varepsilon : \R^d \rightarrow \R$ with radius $\varepsilon > 0$ (see for example \cite[Chapter 4.4]{Brezis}) we have
\begin{equation*}
    u^\varepsilon := u \ast \varphi_\varepsilon \in \mcal{C}^\infty(\mcal{D}).
\end{equation*}
Under the assumption that 
\begin{equation} \label{example_spde_stay_away_from_zero}
    \int_\mcal{D} |u^\varepsilon(x)| \dif x > 0,
\end{equation}
we define
\begin{equation*}
     \phi(x\cond u) := |u^\varepsilon(x)| \, \Big(\int_\mcal{D} |u^\varepsilon(x)| \dif x \Big)^{-1},
\end{equation*}
and thereby get a probability density on $\markspace$ with corresponding distribution $\Phi(\, \cdot \, |u) = \int_\cdot \phi(x|u) \, \dif x$ for any non trivial $u \in \mcal{H}$. \newline

Given a signal path $X$ according to \eqref{example_spde_porous}, under an analogous assumption to \eqref{example_spde_stay_away_from_zero} we define the observation $Y$ as the marked Cox process with $\Prob$-local characteristics $(\lambda^g(X(t)), \Phi(\dif x \cond  X(t)))$. The ground process $Y^g$ is indeed a Cox process in time, as $\lambda^g(X(\, \cdot \, )$ is continuous and $\mcal{F}_0$-measurable by construction, hence \cite[Theorem 14.6.I.]{DVJ_vol2} applies.
\flushright $\diamond$
\end{example}
\section{The Filtering Equations} \label{section_filtering_equations}
In this section we are going to derive the classical equations of the unnormalized and normalized filters for the observation scheme  \ref{introduction_MPP_observation}. The main techniques for this are known since Snyder's seminal paper \cite{Snyder_filtering} and have been generalized to the MPP case by Brémaud, see \cite{Bremaud}. Other references covering the topic are for example \cite{DVJ_vol2, Lipster_shiryaev}. Our paper is the first to tackle the case of an SPDE signal and thus, in comparison to the rather recent paper \cite{Wei_Sun_filtering_MPP}, we do know the explicit form of the generator $\gen$ and the functional analytical framework of $X$. For the rest of this section we assume that Assumptions \ref{ass_signal} and \ref{ass_observation} hold true.\\
 
 \subsection{The Kallianpur-Striebel formula} \label{subsec_change_of_measure_MPP}
As usual in filtering, our first step is to show the existence of a reference measure $\Qrob$ on $(\Omega, \mcal{F})$ under which the process $Y$ has $(\Qrob, \filtt)$-local characteristics $(1, |\markspace|^{-1} \mu_\markspace (\dif x))$, in other words under which $Y$ has a unit rate Poisson-distributed ground process and  uniformly distributed marks in $\markspace$. By $\Prob_t$ and $\Qrob_t$ we denote the restrictions of the respective measures to $\filtt$, for any $t \in [0,T]$. \\ 

First we define the process $\tgeq{\hat{Z}(t)}$ via
\begin{align} \label{radon_nikodym_definition}
      \hat{Z}(t) := \exp & \left\{- \int_0^t \hspace{-4pt} \int_\mathcal{K} \log \{ \lambda(s- , x \cond X(s)) \} \right. Y(\dif s, \dif x)  \\
                    &\quad + \left. \int_0^t \hspace{-4pt} \int_\mathcal{K} \big( \lambda(s, x \cond X(s)) -1 \big) \, \mu_\mathcal{K}(\dif x) \dif s \right\},  \quad t \in [0,T],\nonumber
\end{align}
which is well-defined as $\lambda$ is strictly positive and measurable. It can be easily seen that $\hat{Z}$ is stochastic exponential and follows the integral equation
\begin{align} \label{radon_nikodym_equation_1}
      \hat{Z}(t) =  1 + \int_0^t \hspace{-4pt} \int_\markspace \hat{Z}(s-) &( \lambda(s- , x \cond X(s))^{-1} -1) \,  \times \\
      &\times (Y(\dif s, \dif x) -  \lambda(s, x \cond X(s)) \mu_\markspace(\dif x) \dif s ), \nonumber
\end{align}
which can be found with an application of It\^{o}'s formula. The following result is crucial for the filtering equations:
\begin{lemma} \label{lemma_radon_nikodym_martingale}
    The process $\hat{Z}$ given by \eqref{radon_nikodym_definition} is a $(\Prob, \filtt)$-martingale.
        \begin{flushright}
    $\diamond$
\end{flushright}
\end{lemma}
We omit a detailed proof, as it is standard and widely available in the literature, see 
\cite{Bremaud, DVJ_vol2, Lipster_shiryaev}. Furthermore, the proof does not hinge on the 
specifics of the underlying signal. The general strategy relies on the fact that, as a consequence 
of the boundedness of $Y$, $\hat{Z}$ is a local $(\Prob, \filtt)$-martingale, and by nonnegativity 
also a $(\Prob, \filtt)$-supermartingale. In conclusion, it suffices to show that 
$\E_\Prob [\hat{Z}(t)] = 1$ for any $t \in [0,T]$, under the conditions outlined in 
\cite[VIII.T11]{Bremaud}, which are fulfilled in our case. \\

This lets us introduce the reference probability measure $\dif \Qrob_t := \hat{Z}(t)\dif \Prob_t$, 
which can be extended to a probability measure $\Qrob$ on $(\Omega, \mcal{F})$ by standard methods. 
Under $\Qrob$, the processes $X$ and $Y$ are independent as $Y$ has $(\Qrob, \filtt)$-local 
characteristics $(1, |\markspace|^{-1} \mu_\markspace (\dif x))$, see \cite[VIII.T10]{Bremaud} and  
\cite[Prop. 14.4.III]{DVJ_vol2}. Furthermore, the notion of Radon-Nikodym derivatives is justified and 
we define $\frac{\dif \Qrob_t}{\dif \Prob_t} := \hat{Z}(t)$.

Moreover, as $\hat{Z}$ is nonnegative, we can define $Z(t) = (\hat{Z}(t))^{-1}$, $t \in [0,T]$ , and by \eqref{radon_nikodym_equation_1} get the associated integral equation
\begin{equation} \label{radon_nikodym_equation_2}
      Z(t) =  1 + \int_0^t \hspace{-4pt} \int_\markspace Z(s-) ( \lambda(s- , x \cond X(s)) -1) \,  ( Y(\dif s, \dif x) -   \mu_\markspace(\dif x) \dif s ),
\end{equation} 
for $t \in [0,T].$
Furthermore, the above results imply $\E_\Qrob [Z(t)] = 1$, $t \in [0,T]$ and that the converse Radon-Nikodym derivative is given by $\frac{\dif \Prob_t}{\dif \Qrob_t} := Z(t)$.  \\

\begin{remark}
Depending on the range of the values of \(\lambda\), the canonical choice of \((\Qrob, \filtt)\)-local characteristics \((1, |\markspace|^{-1} \mu_\markspace (\dif x))\) can be adjusted to \((c_g, |\markspace|^{-1} \mu_\markspace (\dif x))\) for some \(c_g > 0\), without any limitations to the theory developed in this paper. All objects derived in this and the subsequent sections can be configured to hold with respect to the adjusted characteristics.

From a numerical perspective, it might be useful to choose \(c_g\) in such a way that the difference \( \lambda(s-, x \cond X(s)) - c_g \) remains within a numerically feasible range in \eqref{radon_nikodym_definition} and forthcoming analogous Radon-Nikodym densities.

From a statistical standpoint, it could be beneficial to choose a \(c_g\) much larger than the actual intensity, analogous to using a reference process with a much higher expected number of points and interpreting the actual observation as a \textit{thinned} point process.
        \begin{flushright}
    $\diamond$
\end{flushright}
\end{remark}

For any  bounded function $\psi: \mcal{H} \rightarrow \R$ we define the \textit{normalized filter} $\tgeq{\eta_t(\psi)}$ by
\begin{equation*}
    \eta_t(\psi) := \E_\mathbb{P} [ \psi(X(t)) | \Yfiltt ], \quad t\in [0,T],
\end{equation*}
where $\tgeq{\Yfiltt}$ is the filtration generated by the observation process $Y$. The starting point of deriving an explicit form for $\tgeq{\eta_t(\psi)}$ is the following Bayes' type formula.

\begin{theorem} \label{lemma_kallianpur_striebel}
The following \textit{Kallianpur-Striebel formula} holds for any  bounded function $\psi: \mcal{H} \rightarrow \R$:
    \begin{equation} \label{kallianpur_striebel_equation}
        \eta_t(\psi) = \frac{\E_\mathbb{Q} [ \psi(X(t)) Z(t) | \mathcal{Y}_t ]}{\E_\mathbb{Q} [  Z(t) | \mathcal{Y}_t ]} \quad \mathbb{P}\text{-a.s.}, \quad t \in [0,T],
    \end{equation}
    where $Z(t)$ is given by
    \begin{align*}
      Z(t) = \exp & \left\{ \int_0^t \hspace{-4pt} \int_\mathcal{K} \log \{ \lambda(s- , x \cond X(s)) \} \right. Y(\dif s, \dif x)  \\
                    &\quad - \left. \int_0^t \hspace{-4pt} \int_\mathcal{K} (\lambda(s, x \cond X(s)) -1) \, \mu_\mathcal{K}(\dif x) \dif s \right\},  \quad t \in [0,T].\nonumber
    \end{align*}
    \begin{flushright}
    $\diamond$
\end{flushright}

\end{theorem}

 \noindent  \textit{Proof.} For any test set $U \in \Yfiltt$ we have for globally bounded $\psi$
\begin{align}
    \E_\Qrob \left[ \1_U \, \E_\Qrob [ \psi(X(t)) Z(t) \cond \Yfiltt]\right] =  \E_\Qrob \left[ \1_U \, \psi(X(t)) Z(t) \right] = \E_\Prob \left[ \1_U \, \psi(X(t)) \right]
\end{align}
by definition and
\begin{align}
        \E_\Qrob \left[ \1_U\, \E_\Prob [  \psi(X(t)) \cond \Yfiltt] \, \E_\Qrob [ Z(t) \cond \Yfiltt]\right] &=  \E_\Prob \left[ \1_U\, \E_\Prob [  \psi(X(t)) \cond \Yfiltt] \right] \\
        &=  \E_\Prob \left[ \1_U\,  \psi(X(t))  \right].
\end{align}
In order to get the equality in ratio form, we observe that for any set $\Yfiltt$-measurable set $N$ on which  $\E_\Prob[Z(t)\cond \Yfiltt] = 0$ we have
\begin{align}
    \Qrob(N) = \E_\Prob[\1_N Z(t)] = \E_\Prob[\1_N \E_\Prob[Z(t)\cond \Yfiltt]] = 0,
\end{align}
implying that \eqref{kallianpur_striebel_equation} holds true under $\Prob$.
The statement for general $\psi$ follows with monotone-class arguments and approximations.
\begin{flushright}
    $\square$
\end{flushright}

\begin{remark} \label{remark_regular_conditional_expectation} To ensure clarity in the notation for regular conditional expectations used in subsequent sections, we define the functional \( z: [0,T] \times \mcal{C}([0,T]; \mcal{H}) \times \mcal{N}^{\# g}_{[0,T] \times \markspace} \rightarrow \R \) by
\begin{align*}
z(t; \mathbf{x}, \xi) := &\exp \Big\{ \int_0^t \hspace{-4pt} \int_\markspace \log \left( \lambda(s,x \mid \mathbf{x}(s)) \right) \xi(\dif s, \dif x) \\
&\qquad \qquad - \int_0^t \hspace{-4pt} \int_\mathcal{K} \left( \lambda(s, x \mid X(s)) - 1 \right) \mu_\mathcal{K}(\dif x) \dif s \Big\},    
\end{align*}
for \( t \in [0,T] \), \( \mathbf{x} \in \mcal{C}([0,T]; \mcal{H}) \), and \( \xi \in \mcal{N}^{\# g}_{[0,T] \times \markspace} \). Given the signal \( X \) and observation \( Y \), we have
\[
Z(t) = z(t; X, Y).
\]
Consequently, the unnormalized posterior distribution is given by
\[
\rho_t(A) := \rho_t(\1_A) = \E_\Qrob [\1_A(X(t)) z(t; X, Y) ] , \quad A \in \Borel(\mcal{H}).
\]
This gives rise to the definition \(\tilde{\rho}: [0,T] \times \mcal{N}^{\# g}_{[0,T] \times \markspace} \rightarrow \mcal{M}^+_{\mcal{H}}\) as follows:
\[
\tilde{\rho}_t \{\chi\} (A) := \E_{X} [ \1_A(X(t)) z(t; \cdot, \xi)] = \E_\Prob[ \1_A(X(t)) \cond Y_{0:t} = \xi_{0:t} ], \quad A \in \Borel(\mcal{H}),
\]
for any \( \xi \in \mcal{N}^{\# g}_{[0,T] \times \markspace} \) and where $\E_{X}$ denotes the expectation under the distribution with respect to the law $\Prob_X$ of $X$.  In other words $\tilde{\rho}_t$ is a regular version of the unnormalized conditional expectation $\rho_t$.  Therefore, for a typical observation \( Y \), we have \( \tilde{\rho}_t \{Y\} (A) = \rho_t(A) \), \( t \in [0,T] \).
    \begin{flushright}
    $\diamond$
\end{flushright}
\end{remark}

\subsection{The Zakai equation} \label{subsection_unnormalized_MPP}
As usual in Bayesian estimation theory, we denote the numerator of \eqref{kallianpur_striebel_equation} as 
\begin{equation*}
    \rho_t(\psi) := \E_\mathbb{Q} [ \psi(X(t)) Z(t) | \mathcal{Y}_t ], \quad t \in [0,T],
\end{equation*}
and call the process $\tgeq{\rho_t(\psi)}$ the \textit{unnormalized filter}.
We have the following theorem for the associated filtering equation:
\begin{theorem}[Zakai equation] \label{Zakai_thm}
    For any It\^{o}-function $\psi$ the following equation for the unnormalized filter  holds  
    \begin{align} \label{Zakai_eq_MPP}
    \rho_t(\psi) &= \rho_0(\psi) + \int_0^t \rho_s(\mathcal{L}\psi) \dif s \\
&\quad + \int_0^t \hspace{-4pt} \int_\mathcal{K}  \rho_{s-}((\lambda(s- ,x\cond  \cdot\, ) -1 )\psi) (Y( \dif s , \dif x) - \mu_\mathcal{K}(\dif x) \dif s), \quad \Qrob\text{-a.s.}, \,  \nonumber
\end{align} 
for any  $t\in [0,T]$, where $\gen$ is given by \eqref{signal_generator}.
    \begin{flushright}
    $\diamond$
\end{flushright}
\end{theorem}
\noindent  \textit{Proof.} Let $\psi$ be a globally bounded It\^{o} function. 
For $t \in [0,T]$ we have by It\^{o}'s lemma for variational solutions of SPDE (see \cite[Thm. 1.2]{Pardoux_SPDE_filter}) and by \eqref{radon_nikodym_equation_2} that
\begin{align}
    \psi&(X(t)) Z(t) = \psi(X(0)) + \int_0^t Z(s)  \, _{\mcal{V}^\ast}\langle A(X(t) ) , \text{D}^1 \psi (X(t)  \rangle_\mcal{V} \, \dif s \\
    &\quad +   \int_0^t Z(s) \, \tr \{ \text{D}^2 \psi(X(t) \; (B(X(s))Q^\frac{1}{2})(B(X(s))Q^\frac{1}{2})^\ast \, \} \, \dif s \nonumber \\
    &\quad +  \int_0^t ( Z(s)D^1\psi(X(s)),  B(X(s)) \dif W(s) )_\mcal{H} \nonumber \\
    &\quad + \int_0^t \hspace{-4pt} \int_\markspace Z(s-) ( \lambda(s- , x \cond X(s)) -1)  \psi(X(s)) \,  [ Y(\dif s, \dif x) -   \mu_\markspace(\dif x) \dif s ]. \nonumber
\end{align}
We take conditional expectations w.r.t. $\Yfiltt$ on both sides and use the definition of the infinitesimal generator in \eqref{signal_generator} to arrive at
\begin{align}
       \E_\Qrob [\psi(X(t)) Z(t)\cond \Yfiltt] &= \E_\Qrob \left[ \psi(X(0))\cond \Yfiltt\right]   \\
       & \quad +  \E_\Qrob \left[ \int_0^t Z(s) \gen(\phi(X(s)) ) \, \dif s \cond \Yfiltt \right]  \nonumber \\
    &\quad + \E_\Qrob \Big[ \int_0^t \hspace{-4pt} \int_\markspace Z(s-) ( \lambda(s- , x \cond X(s)) -1)  \psi(X(s)) \, \times \nonumber \\ 
    &\qquad \qquad  \qquad \qquad \qquad  \quad  \times ( Y(\dif s, \dif x) -   \mu_\markspace(\dif x) \dif s ) \cond \Yfiltt \Big], \nonumber
\end{align}
as the stochastic integral vanishes due to being a local $\Qrob$-martingale. Applying the standard stochastic Fubini argument and then inserting the definition of $\rho_t(\psi)$ finishes the proof for globally bounded $\psi$. Using monotone class arguments and approximations, the assertion for a general $\psi$ can be established.
\begin{flushright}
    $\square$
\end{flushright}

\subsection{The Kushner-Stratonovich equation } \label{subsection_normalized_MPP}
Now that we have proven Zakai's equation for the unnormalized filter $\tinT{\rho_t(\psi)}$ in our setting, we can derive an equivalent equation for the normalized filter $\tinT{\eta_t(\psi)}$ from \eqref{kallianpur_striebel_equation}.  

\begin{corollary}[Kushner-Stratonovich equation] \label{kushner_strat_thm_mpp}
    For any It\^{o}-function $\psi$ the following equation for the normalized filter holds  
    \begin{align} \label{kushner_strat_eq_mpp}
            \eta_t(\psi) &= \eta_0(\psi) + \int_0^t \eta_s(\mathcal{L}\psi) \dif s \\
                    &\quad + \int_0^t \int_\mathcal{K} \frac{\eta_{s-}(\psi \, \lambda(s-, x \cond \cdot))- \eta_{s-}(\psi)\eta_{s-}(\lambda(s-, x \cond \cdot))}{\eta_{s-}(\lambda(s-, x \cond \cdot))} \times \nonumber \\ 
                    & \qquad \qquad \qquad \qquad  \times (Y(\dif s \times \dif x) - \eta_{s-}(\lambda(s-, x \cond \cdot)) \mu_\mathcal{K}(\dif x) \dif s), \nonumber
    \end{align}
    for any $t \in [0,T]$, where $\gen$ is given by \eqref{signal_generator}.
    \begin{flushright}
    $\diamond$
\end{flushright}
\end{corollary}

\begin{proof}
 Let $\psi$ be a globally bounded It\^{o} function.
As usual in filtering theory, we are going to use
\begin{equation} \label{kushner_strat_pf_equation_unnormalized_fraction}
    \eta_t(\psi) = \frac{\rho_t(\psi)}{\rho_t(\1)}, \quad t \in [0,T].
\end{equation}

\noindent As $Z(t)^{-1} =  \hat{Z}(t)$, by \eqref{radon_nikodym_equation_1} we have
\begin{align} \label{kushner_strat_pf_denominator_time_evolution}
    (Z(t))^{-1} &= 1 - \int_0^t \int_\markspace \frac{\lambda(s-, x \cond X(s)) -1}{Z(s-) \, \lambda(s-, x \cond X(s))} \times \\ 
    &\qquad \qquad \times (Y(\dif s, \dif x) -  \lambda(s, x \cond X(s)) \mu_\markspace(\dif x) \dif s ). \nonumber
\end{align}
\noindent From here it can be easily derived that the denominator in \eqref{kushner_strat_pf_equation_unnormalized_fraction} suffices
\begin{align}
    \dif \rho_t(\1)^{-1} &= -  \int_\markspace \frac{\eta_{t-}(\lambda(t-, x \cond \cdot )) -1}{\rho_{t-}(\1) \, \eta_{t-}(\lambda(t-, x \cond \cdot ))} \times \\ 
    &\qquad \qquad \times (Y(\dif t, \dif x) -  \eta_{t}(\lambda(t, x \cond \cdot )) \mu_\markspace(\dif x) \dif t), \nonumber 
\end{align}
see e.g. \cite{Bremaud, DVJ_vol2} for detailed discussions on restrictions of stochastic intensities to smaller filtrations.
Now, an application of It\^{o}'s lemma yields
\begin{align}
    \dif( \rho_t(\psi) &\rho_t(\1)^{-1})  \nonumber\\ 
    &=   \rho_{t-}(\1)^{-1} \dif \rho_t(\psi) + \rho_{t-}(\psi) \dif \rho_{t}(\1)^{-1} + \Delta \rho_t(\psi) \Delta \rho_t(\1)^{-1} \\ 
    &=   \rho_{t}(\1)^{-1}  \rho_t(\mathcal{L}\psi) \dif t \nonumber \\
    &\quad + \int_\markspace  \rho_{t-}(\1)^{-1}  \rho_{t-}((\lambda(t- ,x\cond  \cdot\ ) -1)\psi) (Y( \dif t , \dif x) - \mu_\mathcal{K}(\dif x) \dif t) \nonumber \\ 
    &\quad - \int_\markspace  \rho_{t-}(\psi) \frac{\eta_{t-}(\lambda(t-, x \cond \cdot )) -1}{\rho_{t-}(\1) \, \eta_{t-}(\lambda(t-, x \cond \cdot ))} \times \nonumber \\ 
    &\qquad \qquad \qquad \qquad \qquad \quad \times (Y(\dif t, \dif x) -  \eta_{t}(\lambda(t, x \cond \cdot )) \mu_\markspace(\dif x) \dif t) \nonumber \\
    &\quad - \int_\markspace   \frac{\eta_{t-}(\lambda(t-, x \cond \cdot )) -1}{\rho_{t-}(\1) \, \eta_{t-}(\lambda(t-, x \cond \cdot ))} \rho_{t-}((\lambda(t- ,x\cond  \cdot\ ) -1)\psi) Y( \dif t , \dif x),\nonumber
\end{align}
where all terms are well-defined due to our boundedness assumptions on $\lambda$ and $\psi$. Rearranging terms and inserting the equality \eqref{kushner_strat_pf_equation_unnormalized_fraction} lead to \eqref{kushner_strat_eq_mpp}.
The claim for any $\psi$ follows from monotone-class reasoning and approximation methods.
\end{proof}

\begin{remark} \label{remark_uniqueness_filtering_eq}
The discussion of uniqueness for solutions to the Zakai and Kushner-Stratonovich equations is beyond the scope of this paper. However, under Assumptions \ref{ass_signal}, \ref{ass_observation} and Assumption \ref{ass_intensity_cts}(a) which is going to be introduced in the next section, uniqueness can be established using standard arguments based on the associated martingale problem. In particular, similar techniques to those presented in \cite{ceci2000filtering, ceci2014zakai}  which are again based on \cite[Chapter 4]{ethier1986markov} can be employed to rigorously verify uniqueness. 
\end{remark}

\subsection{The filtering equations for finite dimensional observations}\label{subsection_finite_dim_filtering_eq}

Fix $M \in \N \backslash \{0\}$ and let $Y^M$ be the process from the observation scheme 
\ref{introduction_MPP_discretized_observation_multivariate}. As mentioned in Section 
\ref{sec_2},  $Y^M = (Y^M_1, \dots, Y^M_M)$ is a multivariate point process on 
$[0,T]$ with conditional intensities $(\lambda^M_1(t), \dots, \lambda^M_M(t))$ under $\Prob$. By 
construction, any properties which follow from Assumptions \ref{ass_observation} carry over to the 
counterparts for multivariate point processes.

Generally speaking, the theory of filtering for point processes (without marks) is well-established. However, since there is no known literature addressing the filtering of multivariate point processes with SPDE signals, except for the conference paper by Florchinger \cite{Florchinger_filtering_spde}, we present the main results in this section for the sake of completeness. As a notational convention, we will use the superscript \( M \) to distinguish between finite- and infinite-dimensional objects. \newline

For any $i = 1, \dots, M$, we define
\begin{align*}
\hat{Z}^M_i(t) := \exp \Big\{ \hspace{-3pt} - \hspace{-3pt} \int_0^t \log& \left\{\frac{\lambda^M_i(t\cond X(t))}{\mu_\markspace(K^M_i)} \right\} \, \dif Y^M_{i}(t) +  \\
&\qquad \qquad  +  \int_0^t \big( \lambda^M_i(t\cond X(t)) - \mu_\markspace(K^M_i) \big) \, \dif s \Big\},
\end{align*}
and
\begin{equation*}
    \hat{Z}^M (t) := \prod_{i = 1}^M  \hat{Z}^M_i(t)
\end{equation*}
for any $t \in [0,T]$.\newline

\noindent Analogous to Lemma \ref{lemma_radon_nikodym_martingale}, we have
\begin{lemma} \label{lemma_Z_is_martingale_multivariate}
    The process $\tinT{\hat{Z}^M (t)}$ is a $(\Prob, \filtt)$-martingale.
        \begin{flushright}
    $\diamond$
\end{flushright}
\end{lemma}
As this can be shown by standard techniques, we again omit the proof and  refer to \cite{Bremaud, Lipster_shiryaev, DVJ_vol2}. 

Using above Lemma, analogously to Section \ref{subsec_change_of_measure_MPP} we define the reference probability measure $\dif \Qrob_t^M := \hat{Z}^M (t) \dif \Prob$, which can be extended to a probability measure $\Qrob^M$ on $(\Omega, \mcal{F})$. Under $\Qrob^M$ the process $\tinT{Y^M(t)}$ is an $M$-dimensional Poisson process on $[0,T]$ with rate $\mu_\markspace(K^M_i)$ independent of $X$.

Furthermore, Lemma \ref{lemma_Z_is_martingale_multivariate} implies the existence of the reverse Radon-Nikodym-derivative $\tinT{Z^M (t)}$ by setting $Z^M(t):=(\hat{Z}^M(t))^{-1}$ and that $\tinT{Z^M (t)}$ is a $(\Qrob^M, \filtt)$-martingale as $\E_{\Qrob^M} [\hat{Z}^M(t)] = 1$ for any $t \in [0,T]$. 

Denote by $\tinT{\Yfiltt^M}$ the filtration generated by $\tinT{Y^M(t)}$. We have

\begin{theorem} \label{lemma_kallianpur_striebel_multivariate}
The following \textit{Kallianpur-Striebel formula} holds $\Prob$-a.s. for any  bounded function $\psi: \mcal{H} \rightarrow \R$:
    \begin{equation} 
        \eta^M_t(\psi) := \E_\Prob [\psi(X(t))\cond \Yfiltt^M] = \frac{\E_{\Qrob^M} [ \psi(X(t)) Z^M(t) | \mathcal{Y}^M_t ]}{\E_{\Qrob^M} [  Z^M(t) | \mathcal{Y}^M_t ]}, \quad t \in [0,T],
    \end{equation}
    where $Z^M$ is given by
    \begin{align*}
       Z^M(t) := \exp \Big\{ \sum\limits_{i=1}^M \Big[  &\int_0^t \log \left\{\frac{\lambda^M_i(s\cond X(s))}{\mu_\markspace(K^M_i)} \right\} \, \dif Y^M_{i}(s)  \\       &-  \int_0^t \big(\lambda^M_i(s\cond X(s)) - \mu_\markspace(K^M_i)\big) \, \dif s \Big] \Big\}, \quad t \in [0,T].
    \end{align*}
    \begin{flushright}
    $\diamond$
\end{flushright}
\end{theorem}
\noindent The proof works exactly as the on for Lemma \ref{lemma_kallianpur_striebel} after replacing $Z$, $\Qrob$ and $\Yfiltt$ with their corresponding counterparts with superscript $M$.

We define $\rho^M_t(\psi) := \E_\mathbb{Q^M} [ \psi(X(t)) Z^M(t) | \mathcal{Y}^M_t ]$  and have the analogous results:
\begin{theorem}[Zakai equation for multivariate point processes] \label{Zakai_thm_multivariate}
    For any It\^{o}-function $\psi$  the following equation for the unnormalized filter  holds  
    \begin{align} \label{Zakai_eq_multivariate}
    \rho^M_t(\psi) &= \rho^M_0(\psi) + \int_0^t \rho^M_s(\mathcal{L}\psi) \dif s \\
&\quad +  \sum_{i=1}^M \int_0^t  \rho^M_{s-}((\lambda^M_i(s- \cond  \cdot \, ) -1 )\psi) (Y^M_i( \dif s ) -\mu_\markspace(K^M_i)\dif s), \quad \Qrob^M\text{-a.s.}, \,  \nonumber
\end{align} 
for any  $t\in [0,T]$, where $\gen$ is given by \eqref{signal_generator}.
    \begin{flushright}
    $\diamond$
\end{flushright}
\end{theorem}

\begin{corollary}[Kushner-Stratonovich equation for multivariate point processes] \label{kushner_strat_thm_multivariate}
    For any It\^{o}-function $\psi$   the following equation for the normalized filter holds $\Prob$-a.s. 
    \begin{align} \label{kushner_strat_eq_multivariate}
            \eta^M_t(\psi) &= \eta^M_0(\psi) + \int_0^t \eta^M_s(\mathcal{L}\psi) \dif s \\
                    &\quad + \sum_{i = 1}^M \int_0^t  \frac{\eta^M_{s-}(\psi \, \lambda^M_i(s- \cond \cdot \, ))- \eta^M_{s-}(\psi)\eta^M_{s-}(\lambda^M_i(s- \cond \cdot \,))}{\eta^M_{s-}(\lambda^M_i(s- \cond \cdot \,))} \times \nonumber \\ 
                    & \qquad \qquad \qquad \qquad  \qquad \times (Y^M_i(\dif s ) - \eta^M_{s-}(\lambda^M_i(s- \cond \cdot \, )) \mu_\markspace(K^M_i) \,  \dif s), \nonumber
    \end{align}
    for any $t \in [0,T]$, where $\gen$ is given by \eqref{signal_generator}.
    \begin{flushright}
    $\diamond$
\end{flushright}
\end{corollary}

Both Theorem \ref{Zakai_thm_multivariate} and Corollary \ref{kushner_strat_thm_multivariate} can be proven analogously to Theorem \ref{Zakai_thm} and Corollary \ref{kushner_strat_thm_mpp} by replacing the MPP objects with their multivariate counterparts. For further details we refer to \cite{Bremaud} and \cite{Lipster_shiryaev}.

%
%
%
%
%
%
%
%
%
%
%
%
%
%

%
\section{Consistency of finite-dimensional approximations and error bounds} \label{section_observation_schemes}
In this section, we explore the relationship between the observations from \ref{introduction_MPP_observation} and \ref{introduction_MPP_discretized_observation_multivariate}, as well as the corresponding estimators for the unnormalized and normalized posterior distributions. If we consider \ref{introduction_MPP_observation} as an observation scheme with an "infinitely high" resolution, and \ref{introduction_MPP_discretized_observation_multivariate} as an approximation with limited spatial information, it naturally raises questions about the error bounds between them. To address these questions, we introduce the concept of \textit{dissecting systems}, which are nested partitions commonly used in measure theory.

Using this framework, we construct a nested series of multivariate observations that can be embedded into the MPPs. We demonstrate that this series weakly converges to the process corresponding to \ref{introduction_MPP_observation} in the space of MPPs.

Additionally, we examine the convergence of the corresponding estimators for the normalized and unnormalized posterior distributions. We establish convergence in total variation and provide error bounds.

In the third subsection, we introduce the concept of partial finite-dimensional observation, motivated by the application to CLSM data, where we never observe the entire spatial area but only a fixed subset of partition sets. We derive error bounds for the unnormalized and normalized posterior distributions given these partial observations. 

\label{subsec_convergence_of_posterior_distributions}
\paragraph{Nested partitioning of the markspace}

In order to investigate convergence properties of a family of observation paths according to \ref{introduction_MPP_discretized_observation_multivariate}, $M \in \N$, we have to make assumptions about the underlying corresponding partitions $\markspace^M$,  introduced in Section \ref{sec_2}. The concept of \textit{dissecting systems}, introduced below, is particularly useful for this purpose. It defines a system of nested partitions that interacts well with point process theory and is intuitive to understand. The following definition is taken from \cite{DVJ_vol2}.
\begin{definition} \label{definition_dissecting_system}
    A sequence $(\markspace^M)_{M \in \N}$ of partitions $\markspace^M = \{ K^M_1, \dots , K^M_{n_M}\} $, $M\in \N$, consisting of sets in $\Borel(\markspace)$, is a \textit{dissecting system} for $\markspace$ iff
    \begin{enumerate}[(i)]
        \item The sets $K^M_1, \dots , K^M_{n_M}$ are disjoint and $\bigcupplus\limits_{i=1}^{n_M} K^M_i = \markspace$ for any $M\in \N$.
        \item The $\markspace^M$ are nested with increasing $M$, i.e. $K^{M-1}_i \cap K^{M}_j = K^{M}_j$ or $\emptyset$.
        \item Given any distinct $x_1, x_2 \in \markspace$, there exists a $\tilde{M}\in \N$, such that $x_1 \in K^{\tilde{M}}_i$ implies $x_2 \notin K^{\tilde{M}}_i$.
    \end{enumerate}
       \begin{flushright}
    $\diamond$
\end{flushright} 
\end{definition}

The last property is called the \textit{point-separation property} of the dissecting system. It implies that for any $x \in \markspace$ there exists a uniquely determined nested sequence of sets $(K^M\{x\})_{M\in \N}$ with
\begin{equation*}
    x\in K^M\{x\} \textnormal{ and } K^M\{x\} \in \markspace^M \textnormal{ for any } M \in \N,
\end{equation*}
such that $\bigcap\limits_{M = 1}^\infty K^M\{x\} = \{x\}$. 

As $(K^M\{x\})_{M\in \N}$ is a monotonic sequence, for any measure $\xi$ on $(\markspace, \Borel(\markspace))$ we get by continuity from above that
\begin{equation} \label{dissecting_system_cty_from_above}
    \xi( K^M\{x\} ) \rightarrow \xi(\{x\}) \text{ for } M \rightarrow \infty.
\end{equation}

The markspace $\markspace$ contains a dissecting system, as any Polish space contains at least one, see \cite[Proposition A2.1.IV.]{DVJ_vol1}. Moreover, as $\markspace$ is compact hence bounded, we naturally have that all the sets inside its dissecting systems are bounded.

A practical interpretation of such a dissecting system is seeing $M$ as a theoretically increasing resolution of an image and $\markspace^M$ as the corresponding collection of pixels.

Let $(\markspace^M)_{M \in \N}$ be a fixed dissecting system of $\markspace$ for the rest of this section. Definition \ref{definition_dissecting_system} implies that we may assume the existence of some strictly decreasing real positive sequence $(D_M)_{M  \in \N}$ converging to zero such that
\begin{equation} \label{sec_4_partition_sets_bounded}
    \overline{\textnormal{diam}}(\markspace^M) := \max\limits_{i \leq n_M} \, 
    \textnormal{diam}_{\markspace}(K^M_i) \leq D_M \quad  \text{ for all } M \in \N,
\end{equation}
where  the separability of the underlying space assures the existence of such a dissecting system and $\markspace$ is equipped with the standard metric on $\R^d$.\newline

\subsection{Convergence of finite-dimensional observations} 

\paragraph{Induced MPP} 
In order to discuss the convergence of measures, we need to specify a common measure space.  For a fixed $M\in \N$ let  $Y^{M}$ be the process from the multivariate observation scheme (O$^{M}$), where the according sets $K^{M}_1, \dots, K^{M}_{n_{M}} \in \markspace^{M}$. In particular,  $Y^{M}$ is an $M$-variate counting measure on $[0,T]$, whereas $Y$ is a measure on $[0,T]\times \markspace$. We will demonstrate how the explicit construction in Section \ref{subsection_observation_structures} induces a marked point process (MPP) on the product measure space.

First,  we observe that for any $M \in \N$ we can choose a set of points 
\begin{equation*}
    \mathbf{k}^M := \{ k^M_1 , \dots, k^M_{{n_M}} \in \markspace \, \cond \, k^M_i \in K^M_i, \; i =1,\dots, n_M\}
\end{equation*}
 which we call \textit{representative points} of the corresponding sets. We assume these representative points are chosen by some deterministic rule and that they lie in the inner of the corresponding sets, e.g. choosing the center of each set. By above assumption on the diameter of the partition sets we have
\begin{equation} \label{sec_4_metric_on_markspace_upper_bound}
    d_\markspace(x,k^M_i) \leq  \textnormal{diam}_{\markspace}(K^M_i) \leq D_M\quad  \text{ for all } x \in K^M_i,
\end{equation}
for any $i = 1, \dots, n_M$. Now, let $(\mathbf{k}^M)_{M\in \N}$ be a fixed sequence of representative points  for $(\markspace^M)_{M \in \N}$. Given $Y^{M}$ from (O$^{M}$), we define the MPP $\Tilde{Y}^{M}$ using the representative points by setting
\begin{equation} \label{construction_approx_MPP}
    \Tilde{Y}^{M}( \dif t , \dif x) := \sum_{i = 1}^{M} \sum_{\tau_i \in Y^{M}_i( [0,T])}  \delta_{\tau_i \times k^{M}_i} (\dif t, \dif x), \quad \Prob\text{-a.s.},
\end{equation}

with $ k^{M}_1, \dots,  k^{M}_{n_{M}} \in \mathbf{k}^{M}$. 

The process $\Tilde{Y}^M$ is a re-embedding of the MPP $Y$ and it can be easily seen that it is indeed 
an MPP according to Def. \ref{definition_point_process}(iii). In particular, we observe that the ground processes coincide, i.e.
\begin{equation} 
    (\Tilde{Y}^M)^g = Y^g.
\end{equation}
One can view \(\Tilde{Y}^M\) as an approximation of \(Y\), where \(\Tilde{Y}^M\) does not capture the exact positions of the marks \(\kappa_i\) but only identifies the partition set \(K^M(\kappa_i)\) in which they lie.\\

\paragraph{Weak convergence of observations}
 In the following we are going to explicitly use that $Y$, and therewith also $\Tilde{Y}^M$, are $\Prob$-a.s. totally bounded which implies that we can use the notion of weak convergence instead of weak-hash convergence for boundedly finite measures. For the sake of deriving explicit convergence rates, following \cite{Dudley_Baire} we introduce the space $\textnormal{BL}([0,T]\times \markspace)$ of all bounded Lipschitz functions on $[0,T]\times \markspace$ with the norm
\begin{equation*}
    \| f \|_{\textnormal{BL}} := \| f \|_L + \|f \|_\infty, \quad f \in \textnormal{BL}([0,T]\times \markspace),
\end{equation*}
where 
\begin{equation*}
    \| f \|_L := \sup \left\{ \frac{|f(s,x) - f(t,y)|}{ d ((s,x), (t,y))} \, \bigg| \,  d ((s,x), (t,y)) \neq 0\right\}  .
\end{equation*}
Furthermore, each bounded finite signed measure $\mu$ on $([0,T]\times \markspace, \, \Borel([0,T]\times \markspace))$ defines an element of the dual space of $\textnormal{BL}([0,T]\times \markspace)$ with the norm 
\begin{equation} \label{sec_4_definition_of_BL_star_norm}
    \| \mu \|^\ast_\textnormal{BL} := \sup\limits_{f \in \textnormal{BL}([0,T]\times \markspace)} \left\{  | \int_{[0,T]\times \markspace} \, f  \, \dif \mu | \; \big| \; \| f \|_{\textnormal{BL}} = 1  \right\}.
\end{equation}
and by \cite[Theorem 12.]{Dudley_Baire} the weak topology in the space $\mcal{M}^+_{[0,T]\times \markspace}$ of all nonnegative totally bounded Borel measures on the product measure space coincides with the topology defined by $\| \cdot \|^\ast_\textnormal{BL}$ and as as a direct implication, the same applies to the corresponding topologies on $\mcal{N}^\#_{[0,T]\times \markspace}$.\newline

The following result shows that after the re-embedding of the multivariate point processes according to \ref{introduction_MPP_discretized_observation_multivariate}, the approximations weakly converge to the underlying MPP $Y$. 

\begin{proposition}
    
\label{lemma_weak_convergence_of_approx}
Let $(\markspace^M)_{M \in \N}$ be a dissecting system and  $(\mathbf{k}^M)_{M\in \N}$ be a sequence of corresponding representative points. Furthermore, let $Y$ be an MPP on $[0,T] \times \markspace$. For any $M$ we define the approximating MPP $\Tilde{Y}^M$ via the explicit construction in \eqref{construction_approx_MPP}. Then $\Tilde{Y}^M \stackrel{w}{\rightarrow} Y$ $\Prob$-a.s. in $\mcal{N}^{\# g}_{[0,T]\times \markspace}$ for $M \rightarrow \infty$. Furthermore, we have the approximation error
\begin{equation}
     \| \Tilde{Y}^M - Y \|^\ast_\textnormal{BL} \leq Y^g([0,T])  \, \overline{\textnormal{diam}}(\markspace^M).
\end{equation}
\begin{flushright}
    $\diamond$
\end{flushright}
\end{proposition} 

\noindent  \textit{Proof.}
As discussed above,  the process $\Tilde{Y}^M \in \mcal{N}^{\# g}_{[0,T]\times \markspace}$ $\Prob$-a.s. for any $M\in \N$, thus it is also $\Prob$-a.s. an element of $\mcal{N}^\#_{[0,T]\times \markspace}$.

Let $f \in \textnormal{BL}([0,T]\times \markspace)$  with Lipschitz constant $L_f$.
Then we have $\Prob$-a.s. that
\begin{align}
    &\Big| \int_{[0,T]\times \markspace} f(t,x) Y(\dif t, \dif x) -  \int_{[0,T]\times \markspace} f(t,x) \Tilde{Y}^M(\dif t, \dif x)\Big | \\
    &= \Big| \sum\limits_{(\tau_i, \kappa_i)\in Y([0,T]\times \markspace)} f(\tau_i, \kappa_i) - \sum\limits_{(\tau_i,  \kappa_i)\in Y([0,T]\times \markspace)} f(\tau_i, k^M\{ \kappa_i\}) \Big| \\
    &\leq   \sum\limits_{(\tau_i, \kappa_i)\in Y([0,T]\times \markspace)} | f(\tau_i, \kappa_i) - f(\tau_i, k^M\{ \kappa_i\}) | \\
    &\leq   \sum\limits_{(\tau_i, \kappa_i)\in Y([0,T]\times \markspace)} L_f \, d_\markspace(\kappa_i, k^M\{ \kappa_i\})) \leq   \sum\limits_{\tau_i\in Y^g([0,T])} \|f\|_\textnormal{BL} \,\overline{\textnormal{diam}}(\markspace^M)\\ 
    &\leq  Y^g([0,T])  \|f\|_{\textnormal{BL}}  \, D_M \longrightarrow 0, \quad M \rightarrow \infty,
\end{align}
where we used \eqref{sec_4_metric_on_markspace_upper_bound} and the dominating sequence $(D_M)_{M\in \N}$.  Thus $\Tilde{Y}^M \stackrel{w}{\rightarrow} Y$ $\Prob$-a.s. in $\mcal{N}_{[0,T]\times \markspace}$ by the Portemanteau theorem \cite[Thm. 13.16(ii)]{Klenke}.

In particular, we have that the limit process  $Y \in \mcal{N}^{\# g}_{[0,T]\times \markspace}$ $\Prob$-a.s. by assumption and hence  $\Tilde{Y}^M \stackrel{w}{\rightarrow} Y$ $\Prob$-a.s.  in $\mcal{N}^{\# g}_{[0,T]\times \markspace}$.

The approximation error follows directly by choosing $f \in \textnormal{BL}([0,T]\times \markspace)$ from the subset of functions in $\textnormal{BL}([0,T]\times \markspace)$ with  $\|f \|_{\textnormal{BL}} = 1$ and taking the supremum as in \eqref{sec_4_definition_of_BL_star_norm}.
\begin{flushright}
    $\square$
\end{flushright}
We want to remind the reader, that  $\mcal{N}^{\# g}_{[0,T]\times \markspace}$ is in general not closed under weak convergence as accumulation points might appear in the limit even if every element of a sequence is an MPP. However, in our particular setting, we know that the limit process $ Y \in \mcal{N}^{\# g}_{[0,T]\times \markspace}$ $\Prob$-a.s.  allowing us to state the weak convergence in $\mcal{N}^{\# g}_{[0,T]\times \markspace}$.

\subsection{Asymptotic consistency of posterior distributions} 
In this subsection, we investigate the limiting behavior of the unnormalized and normalized posterior distributions with increasing spatial resolution of the underlying partition. Using our explicit construction, we are able to show convergence in total variation. Additionally, we prove that the approximation error decreases linearly with respect to the size of the partition sets.

For the next results we denote by $\mcal{M}^+_{\mcal{H}}$ the space of all totally bounded positive measures on $\mcal{H}$ and make the following additional assumption on the stochastic intensity of $Y$.
\begin{assumption}\label{ass_intensity_cts} \begin{enumerate}[(a)]
    \item     In addition to all properties from Assumption \ref{ass_observation}, the stochastic intensity $\lambda$ of $Y$ is a continuous function on $[0,T]\times\markspace \times \mcal{H}$ and for all $u \in \mcal{H}$ the bounds 
    \begin{equation*}
            \lambda_- := \inf\limits_{(t,x)\in [0,T]\times \markspace} \lambda(t,x \cond u), \quad  \lambda_+ := \sup\limits_{(t,x)\in [0,T]\times \markspace} \lambda(t,x \cond u)
    \end{equation*}
    exist,  such that
    \begin{equation}
    0< \lambda_- \leq \lambda(t,x\cond u )\leq \lambda_+ < \infty,
\end{equation}
for all $(t,x) \in [0,T] \times \markspace$ .
\item In addition to (a) let the stochastic intensity $\lambda(\cdot, \cdot \cond X(\cdot)) \in \textnormal{BL}([0,T] \times \markspace)$ $\Prob$-a.s. with  deterministic Lipschitz constant $L_\lambda >0$ such that
    \begin{equation*}
        L_\lambda:= \sup \left\{ \frac{|\lambda(s,x \cond u) - \lambda(t,y \cond u)|}{ d ((s,x), (t,y))} \, \bigg| \,  d ((s,x), (t,y)) \neq 0\right\}  ,
    \end{equation*}
for all $u \in \mcal{H}$.
\end{enumerate}

    \begin{flushright}
    $\diamond$
\end{flushright}
\end{assumption}

Under Assumptions \ref{ass_signal}, \ref{ass_observation}, and \ref{ass_intensity_cts}(a), well-posedness—including strong uniqueness—of each filtering equation follows by standard martingale‐problem techniques. We therefore adopt as a standing assumption that all presented filter equations admit strongly unique solutions. For a detailed treatment of uniqueness via the martingale problem, see \cite{ceci2000filtering, ethier1986markov}.

\paragraph{Asymptotics of unnormalized filter} We now present the main theorem of this section. 
\begin{theorem} \label{thm_Zakai_approximations_strong}
Let Assumptions \ref{ass_signal}, \ref{ass_observation}, \ref{ass_intensity_cts}(a) hold true. Furthermore, let $(\markspace^M)_{M \in \N}$ be a dissecting system and  $(\mathbf{k}^M)_{M\in \N}$ be a sequence of corresponding representative points. Given a signal $X$ according to \textnormal{\ref{introduction_spde_basic_formulation}} let $Y$ be the MPP from observation scheme \textnormal{\ref{introduction_MPP_observation}} and $(Y^{M})_{M\in \N}$ a family of multivariate point processes on $[0,T]$, where each $Y^{M}$ is the process from \ref{introduction_MPP_discretized_observation_multivariate} given $\markspace^M$ and $\mathbf{k}^M$. 

Moreover, let $\rho_t$ and $\rho^M_t$ be the unnormalized posterior distributions from Theorem \ref{Zakai_thm} and Theorem \ref{Zakai_thm_multivariate}, respectively, corresponding to \textnormal{\ref{introduction_MPP_observation}} and \ref{introduction_MPP_discretized_observation_multivariate}, for any $M\in \N$  and $t \in [0,T]$. 
Then, we have the following result:
\begin{enumerate}[(i)]
    \item $\| \rho^M_t - \rho_t \|_{\textnormal{TV}} \longrightarrow 0 $ $\Prob$-a.s. in $\mcal{M}^+_{\mcal{H}}$ for $M \longrightarrow \infty$;
    \item Additionally let Assumption \ref{ass_intensity_cts}(b) hold true. Then, we have the pathwise approximation error
\begin{align} 
   \sup\limits_{t \in (0,T]}  &\| \rho^M_t - \rho_t \|_{\textnormal{TV}} \nonumber \\ &\leq C^1_\rho(T,\lambda,Y)\overline{ \textnormal{diam}}(\markspace^M)  + C^2_\rho(T,\lambda,Y)R(T, \lambda,Y, \overline{ \textnormal{diam}}(\markspace^M)  ) \label{thm_zakai_approx_ii_statement}
\end{align}
with terms defined in \eqref{pf_zakai_pathwise_approx_remainder},\eqref{pf_zakai_pathwise_approx_constant_1} and \eqref{pf_zakai_pathwise_approx_constant_2}.

In particular, $R(T, \lambda,Y, \overline{ \textnormal{diam}}(\markspace^M)  ) = \mcal{O}(\overline{ \textnormal{diam}}(\markspace^M)^2)$ as $M \rightarrow \infty$.
\end{enumerate}
\begin{flushright}
    $\diamond$
\end{flushright}
\end{theorem}
\begin{proof}
Fix $M\in \N$ and $t \in [0,T]$. Moreover, let  $X$ be a signal path and let $Y$ be the MPP given $X$. The statement is clear for $t = 0$, so let $t \in (0,T]$. First, we want to show convergence of the Radon-Nikodym densities. 
For any given continuous path $\mathbf{x} \in \mcal{C}([0,T];\mcal{H})$ we denote by $\mathbf{x_{0:t}}$ the restriction of $\mathbf{x}$ up to time $t$. 
Recall that the Radon-Nikodym density $Z^{M}(t)$ associated to $Y^M$ from Section \ref{subsection_finite_dim_filtering_eq} is given as
\begin{align}
    Z&^{M}(t) = \exp \bigg\{ \sum\limits_{i=1}^{M} \Big[ \int_0^t \log \Big\{ (\mu_\markspace(K^M_i))^{-1} \lambda^{M}_i(s\cond X(s)) \Big\} \,  Y^{M}_i( \dif s)  \Big] \nonumber\\
    &\qquad \qquad \qquad \qquad \qquad \qquad  \qquad \qquad -  \int_0^t \hspace{-4pt} \int_\markspace ( \lambda(s,x\cond X(s)) - 1 ) \,\mu_\markspace(\dif y) \dif s \Big] \bigg\}. \nonumber
\end{align}
To improve readability, for  a given typical path $\mathbf{x} \in \mcal{C}([0,T]; \mcal{H}) $  we define 
 \begin{equation*}
     I(t\cond  \mathbf{x}_{0:t}):= \int_0^t \hspace{-4pt} \int_\markspace ( \lambda(s,x\cond  \mathbf{x}(s)) - 1 ) \,\mu_\markspace(\dif x) \dif s ,\quad t \in [0,T].
 \end{equation*}
and 
 \begin{equation*}
     \theta_{M}(s,x\cond \mathbf{x}(s)) := (\mu_\markspace(K^M\{x \}))^{-1} \int_{K^{M} \{ x \}}\lambda(s, y \cond \mathbf{x}(s))  \, \mu_\markspace(\dif y) .
 \end{equation*}
Analogously to the unnormalized regular conditional expectations in Remark \ref{remark_regular_conditional_expectation}, we define the functionals  $z^M : [0,T]\times \mcal{C}([0,T]; \mcal{H}) \times \mcal{N}^{\# g}_{[0,T]\times \markspace} \rightarrow \R $, $M\in \N$, by
\begin{align*}
    z^M(t; \mathbf{x}, \xi) &:= \exp \bigg\{ \int_0^t \hspace{-4pt} \int_\markspace \log \{ \theta_M(s,x \cond \mathbf{x}(s)) \} \xi(\dif s, \dif x)   -I(t \cond \mathbf{x_{0:t}}) \bigg\},
\end{align*}
for $t \in [0,T]$,  $\mathbf{x} \in \mcal{C}([0,T];\mcal{H})$, $\xi \in \mcal{N}^{\# g}_{[0,T]\times \markspace}$. We can write $Z^{M}(t)$ in terms of the underlying MPP $Y$ as
\begin{align}
        Z^{M}(t) &=  \exp \bigg\{ \int_0^t \hspace{-4pt}  \int_\markspace \log \Big\{ \theta_{M}(s,x\cond X(s)) ) \Big\} \,  Y(\dif s, \dif x)  - I(t\cond  X_{0:t})\bigg\} \label{pf_zakai_approx_restating_ZM} \\ 
        &= z^{M}(t; X, Y) \nonumber,
\end{align}
and note that $I(t\cond  X_{0:t})$ already coincides with the $\mu_\markspace(\dif x) \dif t$-integral in $Z(t)$; compare to the derivation in Section \ref{subsec_change_of_measure_MPP}. By continuity and boundedness of $\lambda$  we have
\begin{equation*}
    \| \theta_M(\cdot ,\, \cdot \cond X(\cdot )) - \lambda(\cdot ,\, \cdot \cond X(\cdot )) \|_\infty \rightarrow 0 \text{ for } M \rightarrow \infty.
\end{equation*}
Hence, as $\lambda( \cdot\, , \cdot \cond X(\cdot))$ is assumed to be uniformly bounded from below away from zero by Assumption \ref{ass_intensity_cts}, we also have $ \log \{\theta_M (\cdot ,\, \cdot \cond X(\cdot )) \} \rightarrow \log\{ \lambda(\cdot,\, \cdot \cond X(\cdot)) $ uniformly for $ M \rightarrow \infty$,  yielding the convergence
\begin{align} \label{pf_conv_unnormalized_convergence_of_radon_proxies_1}
    \int_0^t \hspace{-4pt}  \int_\markspace &\log \{ \theta_M(s,x\cond X(s)) \}  Y(\dif s, \dif x) \\
    &\longrightarrow   \int_0^t \hspace{-4pt} \int_\markspace \log \{ \lambda(s,x\cond X(s)) \} \,  Y(\dif s, \dif x) , \nonumber
\end{align}
and hence also
\begin{align}  \label{pf_conv_unnormalized_convergence_of_radon_proxies_2}
    z^M(t;X,Y) &\longrightarrow z(t;X,Y)
\end{align}
$\Prob$-a.s. and in $L^1$ for $ M \rightarrow \infty$, where $z(t;X,Y)$ was introduced in Remark \ref{remark_regular_conditional_expectation} and represents a functional form of the Radon-Nikodym density $Z(t)$.

Analogously to the definition of $\Tilde{\rho}_t$ in Remark \ref{remark_regular_conditional_expectation}, given the signal path $X$, we define the measure-valued functionals $\tilde{\rho}^M:[0,T] \times \mcal{N}^{\# g}_{[0,T]\times \markspace} \rightarrow  \mcal{M}^+_{\mcal{H}} $, $M \in \N$, by 
\begin{align*}  
    \tilde{\rho}^M_t\{\chi \} (A) &:= \E_{X} [ \1_A(X(t)) \, z^M(t;\, \cdot \,, \chi)] \quad A \in \Borel(\mcal{H}).
\end{align*}
We denote  for the observation $Y^M$ corresponding to  \ref{introduction_MPP_discretized_observation_multivariate} by $\Tilde{Y}^M$ the embedding into the MPPs as explained in \eqref{construction_approx_MPP} and used in the proof of Proposition \ref{lemma_weak_convergence_of_approx}.
We see that
\begin{equation}
    \tilde{\rho}^M_t\{Y \} (A) = \tilde{\rho}^M_t\{ \Tilde{Y}^M \} (A) = \rho^M_t(A),
\end{equation}
hence $\tilde{\rho}^M_t(A)$ is equivalent to a regular version of the unnormalized conditional expectation $\rho^M_t(A)$.

The total variation of  $\tilde{\rho}^M_t\{{\chi} \}$ and $\tilde{\rho}_t\{{\chi} \}$ is given w.r.t. the dominating measure $\Prob_X$, so that we have 
\begin{align} \label{pf_zakai_approx_convergence_of_unnormalized_reg}
    \| \tilde{\rho}^M_t\{{\chi} \} - \tilde{\rho}_t \{\chi \} \|_{\textnormal{TV}} = \dfrac{1}{2} \, \E_X [ | z^M(t;\, \cdot \,, {\chi}) - z(t;\, \cdot\, , \chi) | ] \longrightarrow 0
\end{align}
for $M \rightarrow \infty$ by \eqref{pf_conv_unnormalized_convergence_of_radon_proxies_2}. Hence, for any typical observation path $Y$, we have
\begin{align} \label{pf_thm_zakai_approximation_convergence_final_conv}
    \| \rho_t - \rho^M_t \|_{\textnormal{TV}} \leq \|  \tilde{\rho}^M_t\{Y \}  - \rho_t \|_{\textnormal{TV}} = \|  \tilde{\rho}^M_t\{{Y} \}  - \Tilde{\rho}_t \{ Y \} \|_{\textnormal{TV}}  \longrightarrow 0
\end{align}
$\Prob$-a.s. for $M \rightarrow \infty$, proving (i).\\

For the rest of the proof let us denote $ \lambda(t,x) := \lambda(t,x\cond X(t))$ and $\theta_M(t,x) := \theta_M(t,x\cond X(t))$ for better readability. Furthermore we denote by $h_M:=  \overline{ \textnormal{diam}}(\markspace^M)$ and $Y^g_t := Y^g((0,t])$. 

As $\theta_M(t,x)$ is just the mean value of $\lambda(t, \cdot)$ on the set $K^M \{x\}$, the lower bound $\lambda_-$ naturally applies to $\theta_M(t,x)$ as well. To prove the approximation error in (ii), we first note that for any $(t,x) \in [0,T] \times \markspace$ and $M \in \N$ we have by the mean value theorem and Assumption \ref{ass_intensity_cts}(b)
\begin{align}
    |  \log\{ \theta_M(t,x) \} - \log\{ \lambda(t,x) \} | \leq \frac{|\theta_M(t,x)  - \lambda(t,x)| }{\lambda_-} \leq \frac{L_\lambda h_M}{\lambda_-} =: \tilde{L}_\lambda  h_M,
\end{align}
implying 
\begin{align} \label{pf_zakai_approx_upper_bound_for_Y_integral}
    \Big|   \int_0^t \hspace{-4pt} \int_\markspace \big( \log\{ \theta_M(s,x) \}- \log\{ \lambda(s,x) \} \big) Y(\dif s , \dif x)  \Big|
    \leq Y^g_t\tilde{L}_\lambda  h_M.
\end{align}
To construct an upper bound for $Z(t)$, we define
\begin{align*}
\vartheta_\rho(t,\lambda) := &\max\big\{1, \exp\{- t\, (\lambda_- -1)  \, \mu_\markspace(\markspace)\}  \big\}\\
\lambda_\rho := &\max\{1,\lambda_+\}.    
\end{align*}
In particular,
\begin{equation} \label{pf_zakai_approx_case_distinction_overline_theta}
    \lambda_\rho^{ Y^g_t}  {\vartheta}_\rho(t, \lambda) = 
    \begin{cases}
    \begin{aligned}
          &\exp\{- t\, (\lambda_- -1)  \, \mu_\markspace(\markspace)\} , && 0 < \lambda_- \leq  \lambda_+ \leq 1, \\
        &\lambda_+^{ Y^g_t} \,  \exp\{- t\, (\lambda_- -1)  \, \mu_\markspace(\markspace)\}, && 0 < \lambda_- \leq 1 \leq \lambda_+ < \infty,\\
         &\lambda_+^{ Y^g_t}, && 1 \leq \lambda_- \leq  \lambda_+ < \infty.
    \end{aligned}
    \end{cases}
\end{equation}
where the case distinction for the \(\max\) function depends on the values of the bounds $\lambda_-$ and $\lambda_+$ from Assumption \ref{ass_intensity_cts}.

By definition we have that $\lambda_\rho^{ Y^g_t}{\vartheta}_\rho(t, \lambda, Y)$ is monotonically increasing on \([0,T]\).
From here, we can impose the bound 
\begin{equation} \label{pf_zakai_approx_def_theta_overline}
    Z(t) \leq \lambda_\rho^{ Y^g_t} {\vartheta}_\rho(t, \lambda),
\end{equation}
which is immediately evident from \eqref{radon_nikodym_definition} and the above case distinction. 
Thus, by using above bounds we have
\begin{align}
    &| Z(t) - Z^M(t) | \nonumber \\  
    &= Z(t) \label{pf_zakai_approx_use_identity_in_error} 
     \; \big| 1 - \exp\Big\{   \int_0^t \hspace{-4pt} \int_\markspace \big( \log\{ \theta_M(s,x) \} - \log\{ \lambda(s,x) \} \big) Y(\dif s , \dif x)  \Big\} \big|  \\
    &\leq\lambda_\rho^{ Y^g_t}{\vartheta}_\rho(t, \lambda) \big( \exp\{ Y^g_t\tilde{L}_\lambda  h_M\} -1 \big) .\label{pf_zakai_approx_approx_of_ZM}
\end{align}
A first order Taylor expansion with Lagrange form of the remainder yields
\begin{align*}
     \exp\{ Y^g_t\tilde{L}_\lambda  h_M\} -1 = &Y^g_t\tilde{L}_\lambda  h_M +  \frac{1}{2}\big(\tilde{L}_\lambda  Y^g_t\big)^2 \exp\big\{ \tilde{L}_\lambda  Y^g_t \xi_M \big\} h_M^2\\
     \leq   &Y^g_t\tilde{L}_\lambda  h_M  +  \frac{1}{2}\big(\tilde{L}_\lambda  Y^g_t\big)^2 \exp\big\{ \tilde{L}_\lambda  Y^g_t h_M \big\} h_M^2
\end{align*}
 with $0 < \xi_M < h_M$. Denote by
\begin{equation} \label{pf_zakai_pathwise_approx_remainder}
    R(t, \lambda,Y, h_M) := \exp\big\{ \tilde{L}_\lambda  Y^g_t h_M \big\} h_M^2,
\end{equation}
which satisfies $\mcal{O}(h_M^2)$ as $h_M \rightarrow 0$.
By defining
\begin{align}
    C^1_\rho(t,\lambda,Y):=     &\frac{1}{2}\lambda_\rho^{ Y^g_t} \vartheta_\rho(t,\lambda)  \tilde{L}_\lambda Y^g_t \label{pf_zakai_pathwise_approx_constant_1}\\
    C^2_\rho(t,\lambda,Y):=     &\frac{1}{4}\lambda_\rho^{ Y^g_t} \vartheta_\rho(t,\lambda)  \big(\tilde{L}_\lambda Y^g_t \big) ^2 \label{pf_zakai_pathwise_approx_constant_2}
\end{align}
we obtain
\begin{align} \label{pf_zakai_approx_Z-ZM_final_estimate}
         \frac{1}{2}| Z(t) - Z^M(t) |  \leq  C^1_\rho(t,\lambda,Y)h_M + C^2_\rho(t,\lambda,Y)R(t, \lambda,Y, h_M).
\end{align}
The right hand side  does not depend on $X$, thus we can also bound $\E_X [ | Z^M(t) - Z(t) | ]$ by \eqref{pf_zakai_approx_Z-ZM_final_estimate}. Finally, as $Y^g$, $C^1_\rho(\cdot ,\lambda,Y)$, $C^2_\rho(\cdot,\lambda,Y)$ and $R(t, \lambda,Y, h_M)$ grow monotonically in $t$ and all components are bounded on $[0,T]$ we take the supremum to conclude
\begin{align} 
    \sup\limits_{t\in [0,T]} \| \rho^M_t - \rho_t  \|_{\textnormal{TV}} = &\sup\limits_{t\in [0,T]} \dfrac{1}{2} \, \E_X [ | Z^M(t) - Z(t) | \mid \Yfiltt ] \nonumber\\
    \, \leq &C^1_\rho(T,\lambda,Y)h_M + C^2_\rho(T,\lambda,Y)R(T, \lambda,Y, h_M )\label{pf_zakai_approx_Z-ZM_final_estimate_in_sup},
\end{align}
whereby assertion (ii) is proven.
\end{proof}

As an immediate corollary we provide the averaged error bound under the joint law.

\paragraph{Asymptotics of normalized filter}  As a direct corollary of the preceding results, we obtain the following analogous statement for the normalized posterior distribution. We denote by  $\mcal{P}_{\mcal{H}}$ the space of all probability measures on $\mcal{H}$.
\begin{corollary} \label{corollary_Kushner_strat_approximations_strong}
    Let Assumptions \ref{ass_signal}, \ref{ass_observation}, \ref{ass_intensity_cts}(a) hold true.  Moreover, let $\eta_t$ and $\eta^M_t$ be the unnormalized posterior distributions from Theorem \ref{kushner_strat_thm_mpp} corresponding to \textnormal{\ref{introduction_MPP_observation}} and from Theorem \ref{kushner_strat_thm_multivariate} corresponding to \textnormal{(O$^{M}$)}, respectively, for any $M\in \N$ and $t \in [0,T]$.

Then we have the following result: 
\begin{enumerate}[(i)]
    \item $\| \eta^M_t - \eta_t \|_{\textnormal{TV}} \longrightarrow 0 $  $\Prob$-a.s.  in $\mcal{P}_{\mcal{H}}$ for $M \rightarrow \infty$. 
    \item Additionally let Assumption \ref{ass_intensity_cts}(b) hold true. Then, we have the pathwise approximation error
\begin{align} 
   \sup\limits_{t \in (0,T]}  &\| \eta^M_t - \eta_t \|_{\textnormal{TV}} \nonumber \\ &\leq C^1_\eta(T,\lambda,Y)\overline{ \textnormal{diam}}(\markspace^M)  + C^2_\eta(T,\lambda,Y)R(T, \lambda,Y, \overline{ \textnormal{diam}}(\markspace^M)  ) 
\end{align}
with terms defined in \eqref{pf_zakai_pathwise_approx_remainder},\eqref{pf_kushner_pathwise_approx_constant_1} and \eqref{pf_kushner_pathwise_approx_constant_2}.

In particular, $R(T, \lambda,Y, \overline{ \textnormal{diam}}(\markspace^M)  ) = \mcal{O}(\overline{ \textnormal{diam}}(\markspace^M)^2)$ as $M \rightarrow \infty$.
\end{enumerate}

\begin{flushright}
    $\diamond$
\end{flushright}
\end{corollary}
\begin{proof} Let $\psi \in \mcal{C}_b(\mcal{H})$.
By definition  
\begin{equation*}
    \eta_t(\psi) := \frac{\rho_t(\psi)}{\rho_t(\1)}, \quad \eta^M_t(\psi) := \frac{\rho^M_t(\psi)}{\rho^M_t(\1)}, \; M \in \N,
\end{equation*}
are probability measures on $\mcal{H}$ for all $t \in [0,T]$. By Theorem \ref{thm_Zakai_approximations_strong} we get the convergence $\rho^M_t(\psi)\rightarrow \rho_t(\psi)$ and $\rho^M_t(\1) \rightarrow \rho_t(\1)$ in total variation for  $M \rightarrow \infty$. Moreover,  $\rho^M_t(\1) >0 $  $\Prob$-a.s., thus  $ \eta^M_t(\psi) \rightarrow \eta_t(\psi)$ in total variation by  $\eta_t(\psi)$ being the quotient of two converging sequences. As $\psi$ was chosen arbitrarily assertion (i) follows. \\\\
For the proof of the rate in (ii) we again use the notation of the proof of Thm. \ref{thm_Zakai_approximations_strong} (ii).  As $   \eta_t, \eta^M_t \in \mcal{P}_{\mcal{H}}$, $M\in \N$,
\begin{equation}
    \| \eta^M_t - \eta_t\|_{\textnormal{TV}} = \sup\limits_{A \in \Borel(\mcal{H})} | \eta^M_t(A) - \eta_t(A) |.
\end{equation} 
In order to impose a lower bound we use similar arguments and define the auxiliary constants
\begin{align*}
    \underline{\vartheta}(t, \lambda) := &\min \big\{ 1, \,  \exp\{- t\, (\lambda_+ -1)  \, \mu_\markspace(\markspace)\} \big\},\\
    \underline{\lambda} := &\min \{1, \lambda_-\},
\end{align*}
where an analogous case distinction to \eqref{pf_zakai_approx_case_distinction_overline_theta} is given by
\begin{equation} \label{pf_zakai_approx_case_distinction_underline_theta}
    \underline{\lambda}^{Y^g_t}\underline{\vartheta}(t, \lambda) = 
    \begin{cases}
    \begin{aligned}
          & \lambda_-^{^{Y^g_t}}, && 0 < \lambda_- \leq  \lambda_+ \leq 1, \\
        &\lambda_-^{^{Y^g_t}} \,  \exp\{- t\, (\lambda_+ -1)  \, \mu_\markspace(\markspace)\}, && 0 < \lambda_- \leq 1 \leq \lambda_+ < \infty,\\
         &\exp\{- t\, (\lambda_+ -1)  \, \mu_\markspace(\markspace)\}, && 1 \leq \lambda_- \leq  \lambda_+ < \infty.
    \end{aligned}
    \end{cases}
\end{equation}
It is immediately clear that $\underline{\vartheta}(\cdot, \lambda, Y)$ is monotonically decreasing on $[0,T]$. This gives rise to impose the lower bound
\begin{align}
     Z(t) \geq \underline{\lambda}^{Y^g_t}\underline{\vartheta}(t, \lambda), \label{pf_cor_kushner_approx_underline_vartheta_def}
 \end{align}
for any $ t \in [0,T]$. One can easily verify, that for any $M \in \N$ we also have 
\begin{equation}
   \underline{\lambda}^{Y^g_t}\underline{\vartheta}(t, \lambda) \leq  Z^M(t) \leq {\lambda}_\rho^{Y^g_t}{\vartheta}_\rho(t, \lambda),
\end{equation}
and as the bounds do not depend on $X$ it also follows that
\begin{equation}
         \underline{\lambda}^{Y^g_t}\underline{\vartheta}(t, \lambda) \leq  \rho_t(\1), \rho_t^M(\1) \leq {\lambda}_\rho^{Y^g_t}{\vartheta}_\rho(t, \lambda).
\end{equation}
Hence, for any $A \in \Borel(\mcal{H})$, we conclude
\begin{align}
  | \eta^M_t(A) - \eta_t&(A) | =  \Big| \frac{\rho^M_t(A)}{\rho^M_t(\1)} -  \frac{\rho_t(A)}{\rho_t(\1)} \Big| \nonumber \\
   &=  \Big| \frac{\rho^M_t(A) \, \rho_t(\1) - \rho_t(A) \, \rho^M_t(\1)}{\rho^M_t(\1) \, \rho_t(\1)} \Big| \\
    &\leq (\rho^M_t(\1) \, \rho_t(\1) )^{-1} \, \Big( \,  \big| \, \rho^M_t(A) \, \rho_t(\1) - \rho_t(A) \, \rho_t(\1) \, \big| \\
    &\qquad \qquad \qquad \qquad \qquad  +  \big| \, \rho_t(A) \, \rho^M_t(\1) - \rho_t(A) \, \rho_t(\1) \, \big|  \Big) \nonumber \\
    &\leq \rho^M_t(\1)^{-1}\Big( \,  \big| \, \rho^M_t(A)  - \rho_t(A)  \big|        +  \big|  \rho^M_t(\1) -  \rho_t(\1) \, \big|  \Big)\\
    &  \leq 2\,\Big(\frac{\lambda_\rho}{\underline{\lambda}}\Big)^{Y^g_t}( \underline{\vartheta}^{-1}   {\vartheta}_\rho)(t, \lambda)  \|  \rho^M_t  - \rho_t \|_{\textnormal{TV}} ,
    \label{pf_kushner_approx_TV_bound_in_t}
\end{align}
where we used  that $ \rho_t(\1) \geq  \rho_t(A)$ for all $A  \in \Borel(\mcal{H})$ and the notation $( \underline{\vartheta}^{-1}   {\vartheta}_\rho)(t, \lambda):=  \underline{\vartheta}^{-1}(t, \lambda)   {\vartheta}_\rho(t, \lambda)$ for better readability.  Taking the supremum over all $A \in \Borel(\mcal{H})$ shows that we can bound $\| \eta^M_t - \eta_t \|_{\textnormal{TV}}$ by \eqref{pf_kushner_approx_TV_bound_in_t}.  From here, we estimate using the right hand side of \eqref{pf_zakai_approx_Z-ZM_final_estimate} from the proof of Theorem \ref{thm_Zakai_approximations_strong} to get the stated bound with combined constants
\begin{align}
 &C_\eta^1(t,\lambda,Y):=\lambda_\eta^{Y^g_t} \vartheta_\eta(t,\lambda)\tilde{L}_\lambda Y_t^g   \label{pf_kushner_pathwise_approx_constant_1}, \\
  &C_\eta^2(t,\lambda,Y):=\frac{1}{2}\lambda_\eta^{Y^g_t} \vartheta_\eta(t,\lambda)\big(\tilde{L}_\lambda Y_t^g   \big)^2\label{pf_kushner_pathwise_approx_constant_2}, \\
    &\lambda_\eta := \frac{\lambda_\rho^2}{\underline{\lambda}}, \quad \vartheta_\eta(t,\lambda) := \left( \underline{\vartheta}^{-1}   {\vartheta}^2_\rho\right)(t, \lambda). \nonumber
\end{align}

It can be easilty checked that both constants are monotonically increasing, whence taking the supremum over $t$ finishes the proof of assertion (ii). 
\end{proof}
\begin{remark}
    Exploiting higher Hölder regularity of $\lambda$ in the Kantorovich–Rubinstein distance fails to improve the convergence rate. In fact, the error contribution from the remainder term saturates at the linear rates given in Thm. \ref{thm_Zakai_approximations_strong} and Cor. \ref{corollary_Kushner_strat_approximations_strong}, regardless of further derivatives. Accordingly, we omit the straightforward calculations. Intuitively, this indicates that the accuracy of the filter is governed primarily by the discretization scheme and cannot be substantially enhanced by increased regularity of the underlying intensity.
\end{remark}

\subsection{Partial observations} \label{subsec_partial_observations}

As opposed to our observation schemes models \ref{introduction_MPP_observation} and \ref{introduction_MPP_discretized_observation_multivariate}, where we always have information about the whole mark space $\markspace$, CLSM data does only contain information about a subset of the partition. For the purpose of modeling such a \textit{partial observation scheme}, let $(\markspace^M)_{M \in \N}$ again  be the fixed dissecting system for the mark space $\markspace$ from the last section. For some fixed $M\in \N$, let $\mcal{I}_{M} := \{i_1, \dots i_{|\mcal{I}_{M}|}\}  \subseteq \{1,\dots, n_M\}$ be some subset of indices with $|\mcal{I}_{M}| < n_M$ and let $\markspace^{M}_{\mcal{I}_{M}} $ be the collection of all sets $K^{M}_i$ with $i \in \mcal{I}_{M}$. Because of $|\mcal{I}_{M}| < n_M$ the family $\markspace^{M}_{\mcal{I}_{M}} $ is no longer a partition of $\markspace$.

\paragraph{Partial filtering problem} We can use the tools from Section \ref{subsection_finite_dim_filtering_eq} to derive the analogous filtering equations for a partial observation, as the partition property of $\markspace^M$ is not explicitly required in this context. 

For the signal process $X$ from \textnormal{\ref{introduction_spde_basic_formulation}} we again introduce the $M$-variate observation $Y^{M}$ from \textnormal{(O$^{M}$)} given  $\markspace^{M}$. Now, in addition to that we define the partial observation $Y^{M}|_{\mcal{I}_{M}}$ given the collection of sets $\markspace^{M}_{\mcal{I}_{M}} $. 
Analogously to $Y^M$, we can introduce a  reference measure  $\Qrob^{M}_{\mcal{I}_{M}} $ under which $Y^{M}|_{\mcal{I}_{M}}$ is a $|\mcal{I}_{M}| $-dimensional Poisson process with rate $\mu_\markspace(K^M_i)$ in each component, with Radon-Nikodym derivative $Z^{M}_{\mcal{I}_{M}}(t) := \frac{\dif \Prob|_t}{\dif \Qrob^{M}_{\mcal{I}_{M}}|_t}$ given by
\begin{align*}
Z^{M}_{\mcal{I}_{M}}(t) := \exp &\Big\{ \sum\limits_{i \in \mcal{I}_{M}}  \hspace{-3pt} \int_0^t \log \left\{\frac{\lambda^M_i(s\cond X(s))}{\mu_\markspace(K^M_i)} \right\} \, \dif Y^M_{i}(s) \\
&\qquad \qquad -  \int_0^t \big( \lambda^M_i(s\cond X(s)) - \mu_\markspace(K^M_i) \big) \, \dif s \Big\},
\end{align*}
for any $t \in [0,T]$.
By introducing the filtration $\tinT{\mcal{Y}^{\mcal{I}_{M}}_t}$ generated by  $Y^{M}|_{\mcal{I}_{M}}$ one can derive the unnormalized and normalized posterior distributions $\rho^{\mcal{I}_{M}}_t$ and $\eta^{\mcal{I}_{M}}_t$, respectively, in the exact same way as we did in Section \ref{subsection_finite_dim_filtering_eq}. \newline

The partial observation  $Y^{M}_{\mcal{I}_{M}}$ does not inherit all jumps of $Y$,  only those with marks in the sets of $\markspace^{M}_{\mcal{I}_{M}}$. We were able to interpret the process $Y^{M}$ as an approximation of the MPP $Y$ with uncertainty about the exact mark positions. A crucial property of the embedding $\Tilde{Y}^{M}$  was the identity of the ground processes, i.e. $( \Tilde{Y}^{M})^g = Y^g$, and that  we had $\Tilde{Y}^{M}([0,T]\times \markspace) = Y([0,T]\times \markspace)$. As opposed to that, in general for the partial observation we have  
\begin{equation}
    Y^{M}_{\mcal{I}_{M}} ([0,T]) \leq  Y([0,T]\times \markspace),
\end{equation}
meaning that we may always miss some points.\newline

\paragraph{Approximation error bounds} Although we cannot expect convergence of the estimators in general, we can still derive approximation errors for the total variation distances $\| \rho_t - \rho^{\mcal{I}_{M}}_t \|_{\textnormal{TV}}$ and $\| \eta_t - \eta^{\mcal{I}_{M}}_t \|_{\textnormal{TV}}$, as demonstrated in the next theorem. As we trivially have
\begin{equation} \label{partial_observations_triangle_total_variation}
    \| \rho_t - \rho^{\mcal{I}_{M}}_t \|_{\textnormal{TV}} \leq \| \rho_t - \rho^M_t \|_{\textnormal{TV}} + \| \rho^M_t  - \rho^{\mcal{I}_{M}}_t \|_{\textnormal{TV}},
\end{equation}
for any $t \in [0,T]$, and the analogous inequality for $\| \eta_t - \eta^{\mcal{I}_{M}}_t \|_{\textnormal{TV}} $, the first terms on the right hand sides of the bounds in (i) and (ii) follow direcly by Theorem \ref{thm_Zakai_approximations_strong} and Corollary \ref{corollary_Kushner_strat_approximations_strong}, respectively. Hence, the errors comprise two components: the discretization errors $\kappa_\rho$ from Theorem \ref{thm_Zakai_approximations_strong} and $\kappa_\eta$ from Corollary \ref{corollary_Kushner_strat_approximations_strong}, and additional errors $\epsilon_\rho$ and $\epsilon_\eta$, respectively, that exponentially depend on the size of the unobserved area. The latter accounts for the information loss due to observing only a subset of the partition. We only state the pathwise bounds as the averaged (risk free)w bounds are again a simple corollary of the following Proposition \ref{proposition_error_bounds_partial_observations}.

For better readability we define for any index set $\mcal{I}_M \subset \{1, \dots, M\} $  
\begin{align*}
                   \mcal{I}_M^\tcomplement &:= \{1, \dots, M\} \backslash \mcal{I}_M , \\
                \markspace^M(\mcal{I}_M) &:= \bigcupplus_{i \in \mcal{I}_M} K^M_i .
\end{align*}

\begin{proposition} \label{proposition_error_bounds_partial_observations}
Let Assumptions \ref{ass_signal}, \ref{ass_observation}, \ref{ass_intensity_cts} hold true. Then we have the following pathwise approximation errors.
\begin{enumerate}[(i)]
    \item We have 
    \begin{align*}
        \sup_{t\in [0,T]} \| \rho_t - &\rho^{\mcal{I}_{M}}_t \|_{\textnormal{TV}} \\ \leq &C^1_\rho(T,\lambda,Y)\overline{ \textnormal{diam}}(\markspace^M)  + C^2_\rho(T,\lambda,Y)R(T, \lambda,Y, \overline{ \textnormal{diam}}(\markspace^M)  ) \\
        &+ \epsilon_\rho(\lambda, Y,T)(\markspace(\mcal{I}_M^\tcomplement)),
    \end{align*}
    where $\epsilon_\rho(\lambda, Y,T)$ is defined in \eqref{pf_partial_error_term_unnormalized_form} and the remaining constants are specified in Thm. \ref{thm_Zakai_approximations_strong} (ii).

    \item We have 
    \begin{align*}
            \sup_{t\in [0,T]} \| \eta_t - &\eta^{\mcal{I}_{M}}_t \|_{\textnormal{TV}} \\ \leq &\,C^1_\eta(T,\lambda,Y)\overline{ \textnormal{diam}}(\markspace^M)  + C^2_\eta(T,\lambda,Y)R(T, \lambda,Y, \overline{ \textnormal{diam}}(\markspace^M)  )\\
        &+ \epsilon_\eta(\lambda, Y,T)(\markspace(\mcal{I}_M^\tcomplement)),
    \end{align*}
    where $\epsilon_\eta(\lambda, Y,T)$ is defined in \eqref{pf_partial_error_term_normalized_form} and the remaining terms are provided in Corollary \ref{corollary_Kushner_strat_approximations_strong} (ii) .
\end{enumerate}
\begin{flushright}
    $\diamond$
\end{flushright}
\end{proposition} 
\begin{proof}
Analogously to the proofs of the preceeding approximation error, for (i) we rewrite
\begin{align}
    \| \rho^M_t  &- \rho^{\mcal{I}_{M}}_t \|_{\textnormal{TV}} = \frac{1}{2} \E_X [ | Z^M(t) - Z^{M}_{\mcal{I}_{M}}(t) |  ] \nonumber \\
    &= \frac{1}{2} \E_X \Big[  Z^M(t) \, \big| 1 - \exp\big\{ - \int_0^t \hspace{-4pt} \int_{\markspace^M(\mcal{I}^\tcomplement_M)} \log \{ \theta_M(s,x) \} \, Y(\dif s, \dif x ) \label{pf_partial_zakai_form_of_TV_difference}\\
    &\qquad  \qquad\qquad\qquad\qquad\qquad- \int_0^t \hspace{-4pt} \int_ {\markspace^M(\mcal{I}^\tcomplement_M)} (\lambda(s,x) -1) \, \mu_\markspace(\dif x) \, \dif s   \big\}  \big|  \Big] \nonumber
\end{align}
where $\theta_M(s,x)$ is defined as in the proof of Theorem \ref{thm_Zakai_approximations_strong}. 

As we have $\lambda_- \leq \theta_M(s,x) \leq \lambda_+ $, we can conclude that
\begin{align*}
| \int_0^t \hspace{-4pt} \int_{\markspace^M (\mcal{I}^\tcomplement_M)} &\log \{ \theta_M (s,x) \} \, Y(\dif s, \dif x ) | \\
&\leq \max \big\{  |\log\{\lambda_-\}|, |\log\{\lambda_+\}|   \big\} Y((0,t]\times \markspace^M (\mcal{I}^\tcomplement_M)  ) .
\end{align*}
Hence,
\begin{align}
    &|\int_0^t \hspace{-4pt} \int_{\markspace^M(\mcal{I}^\tcomplement_M)} \log \{ \theta_M(s,x) \} \, Y(\dif s, \dif x ) + \int_0^t \hspace{-4pt} \int_ {\markspace^M(\mcal{I}^\tcomplement_M)} (\lambda(s,x) -1) \, \mu_\markspace(\dif x)| \nonumber \\
        &\quad \leq \max \big\{  |\log\{\lambda_-\}|, |\log\{\lambda_+\}|   \big\} Y((0,t]\times \markspace^M (\mcal{I}^\tcomplement_M)  )  \\ 
        &\qquad \qquad \qquad \qquad  \qquad \qquad + t \mu_\markspace(\markspace^M(\mcal{I}^\tcomplement_M)) \max\{|\lambda_- -1|, |\lambda_+ -1|\} \nonumber
\end{align}
With a similar approximation as in the proof of Theorem \ref{thm_Zakai_approximations_strong}, we now have for given $X$ 
\begin{align}
    &Z^M(t) \, \big| 1 - \exp\big\{ - \int_0^t \hspace{-4pt} \int_{\markspace^M(\mcal{I}^\tcomplement_M)} \log \{ \theta_M(s,x) \} \, Y(\dif s, \dif x ) \nonumber \\
    &\qquad  \qquad\qquad\qquad\qquad\qquad- \int_0^t \hspace{-4pt} \int_ {\markspace^M(\mcal{I}^\tcomplement_M)} (\lambda(s,x) -1) \, \mu_\markspace(\dif x) \, \dif s   \big\}  \big|  \nonumber \\
    &\leq  \lambda_\rho^{Y_t^g} {\vartheta}_\rho(t,\lambda) \, \Big( \exp\Big\{ | \int_0^t \hspace{-4pt} \int_{\markspace^M(\mcal{I}^\tcomplement_M)} \log \{ \theta_M(s,x) \} \, Y(\dif s, \dif x )  \\
    &\qquad  \qquad\qquad\qquad\qquad + \int_0^t \hspace{-4pt} \int_ {\markspace^M(\mcal{I}^\tcomplement_M)} (\lambda(s,x) -1) \, \mu_\markspace(\dif x) \, \dif s |  \Big\}  -1 \Big)  \nonumber \\
    &\leq   \lambda_\rho^{Y_t^g} {\vartheta}_\rho(t,\lambda) \, \Big( \exp\Big\{ \max \big\{  |\log\{\lambda_-\}|, |\log\{\lambda_+\}|   \big\} Y((0,t]\times \markspace^M (\mcal{I}^\tcomplement_M)  )  \label{pf_partial_zakai_final_estimate}\\ 
        &\qquad  \qquad \qquad  \qquad \qquad + t \mu_\markspace(\markspace^M(\mcal{I}^\tcomplement_M)) \max\{|\lambda_- -1|, |\lambda_+ -1|\} \Big\}  -1 \Big)  \nonumber  \\
    &=: \epsilon_\rho(\lambda, Y,t)(\markspace(\mcal{I}_M^\tcomplement))   \label{pf_partial_error_term_unnormalized_form}.     
\end{align}

All components are independent of $X$ and bounded and monotonically increasing on $[0,T]$, hence we can bound \eqref{pf_partial_zakai_form_of_TV_difference} using \eqref{pf_partial_zakai_final_estimate}. Taking the supremum over $t$ proves assertion (i).  \newline

For the proof of the bound in (ii), we observe that
\begin{align*}
    Z^{M}_{\mcal{I}_{M}}(t) \geq  \underline{\lambda }^{Y((0,t] \times \markspace^M(\mcal{I}_M))} 
 \underline{\vartheta}_{\mcal{I}_{M}}(t,\lambda),
\end{align*}
with
\begin{align*}
    \underline{\vartheta}_{\mcal{I}_{M}}(t,\lambda,Y) :=  \min \Big\{ 1 ,  \,  \exp\big\{- t\, (\lambda_+ -1)  \, \mu_\markspace\big(\markspace^M(\mcal{I}_M)\big)\big\}\Big\},
\end{align*}
and where 
\begin{equation*}
    1 \geq \underline{\vartheta}_{\mcal{I}_{M}}(t,\lambda,Y) \geq \underline{\vartheta}(t,\lambda,Y).
\end{equation*}
The remainder of the proof proceeds analogously to the argument presented in the proof of Corollary \ref{corollary_Kushner_strat_approximations_strong} (ii). This yields the following expression for the error term: 
\begin{align}
    \epsilon_\eta(\lambda, Y,t)&(\markspace(\mcal{I}_M^\tcomplement)) \nonumber\\
    &:= \underline{\lambda }^{-Y((0,t] \times \markspace^M(\mcal{I}_M))} 
 \underline{\vartheta}^{-1}_{\mcal{I}_{M}}(t,\lambda) \lambda_\rho^{Y_t^g} {\vartheta}_\rho(t,\lambda) 
 \epsilon_\rho(\lambda, Y,t)(\markspace(\mcal{I}_M^\tcomplement)), \label{pf_partial_error_term_normalized_form}
\end{align}
and, for the sake of brevity, the details are omitted.
\end{proof}
\section{Simulations} \label{sec_simulations}
In this section, we will compare our theoretical results with numerical experiments. The Python code used for the simulations and plots is publicly available at \href{https://github.com/jszala/SPDE_Poisson_filtering.git}{"https://github.com/jszala/SPDE\_Poisson\_filtering.git"}.
\subsection{Synthetic data}
Signal and observation processes both are simulated using explicit Euler schemes in time and finite differences in space. The Git repository also provides the necessary data for reproducing the experiments.

In the experiments, the observation process will be given as a multivariate Poisson process according to the scheme \ref{introduction_MPP_discretized_observation_multivariate}. The intensity is chosen as
\begin{equation} \label{simulations_intensity_form}
    \lambda(t,x) := e^{-at} (c x)^2 \vee \mcal{C}_{\max}
\end{equation}
with $a >0$ being an optional and sufficiently small decay parameter and $c >0$ being a scaling parameter and $\mcal{C}_{\max}$ is some sufficiently large upper bound.  We note that since $\lambda$ is not Lipschitz continuous, the error bounds provided in Corollary \ref{corollary_Kushner_strat_approximations_strong}(ii) are not directly applicable in this case.

\begin{figure}
     \centering
     \begin{subfigure}[b]{1.0\textwidth}
         \centering
         \includegraphics[width=\textwidth]{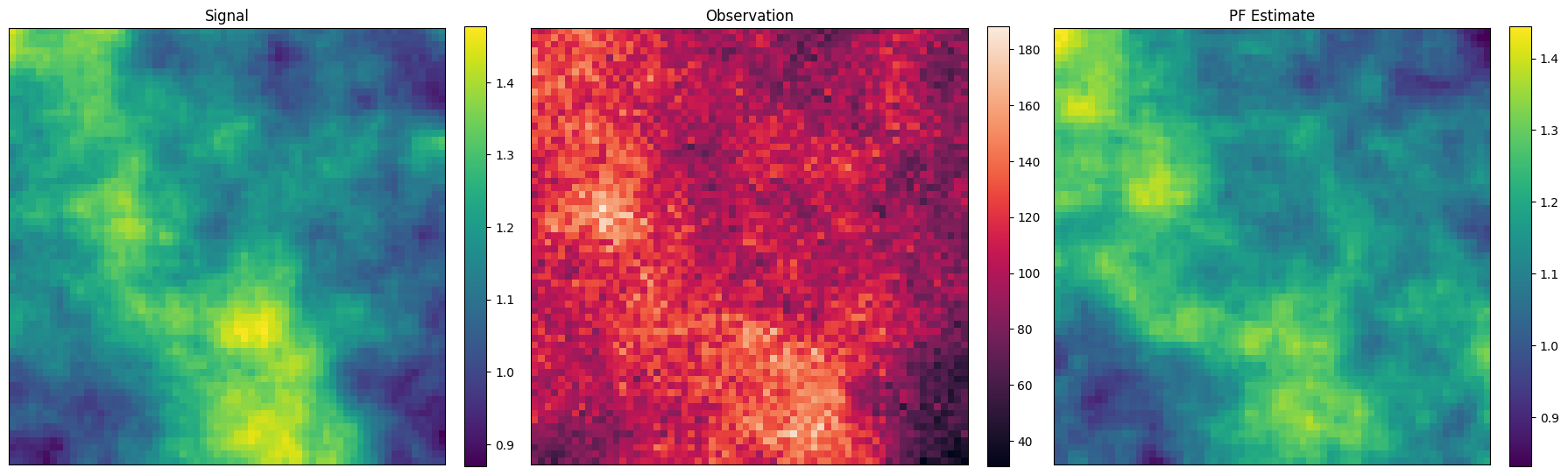}
         \caption{Comparison of state in timestep $t= 192$. True SPDE signal (left), $64\times64$ Poisson observations (center), and particle‐filter posterior mean (right).}
         \label{fig: resolution_figure_SPDE}
     \end{subfigure}
     \hspace{100pt}
     \hfill 
     \begin{subfigure}[b]{1.0\textwidth}
         \centering
         \includegraphics[width=\textwidth]{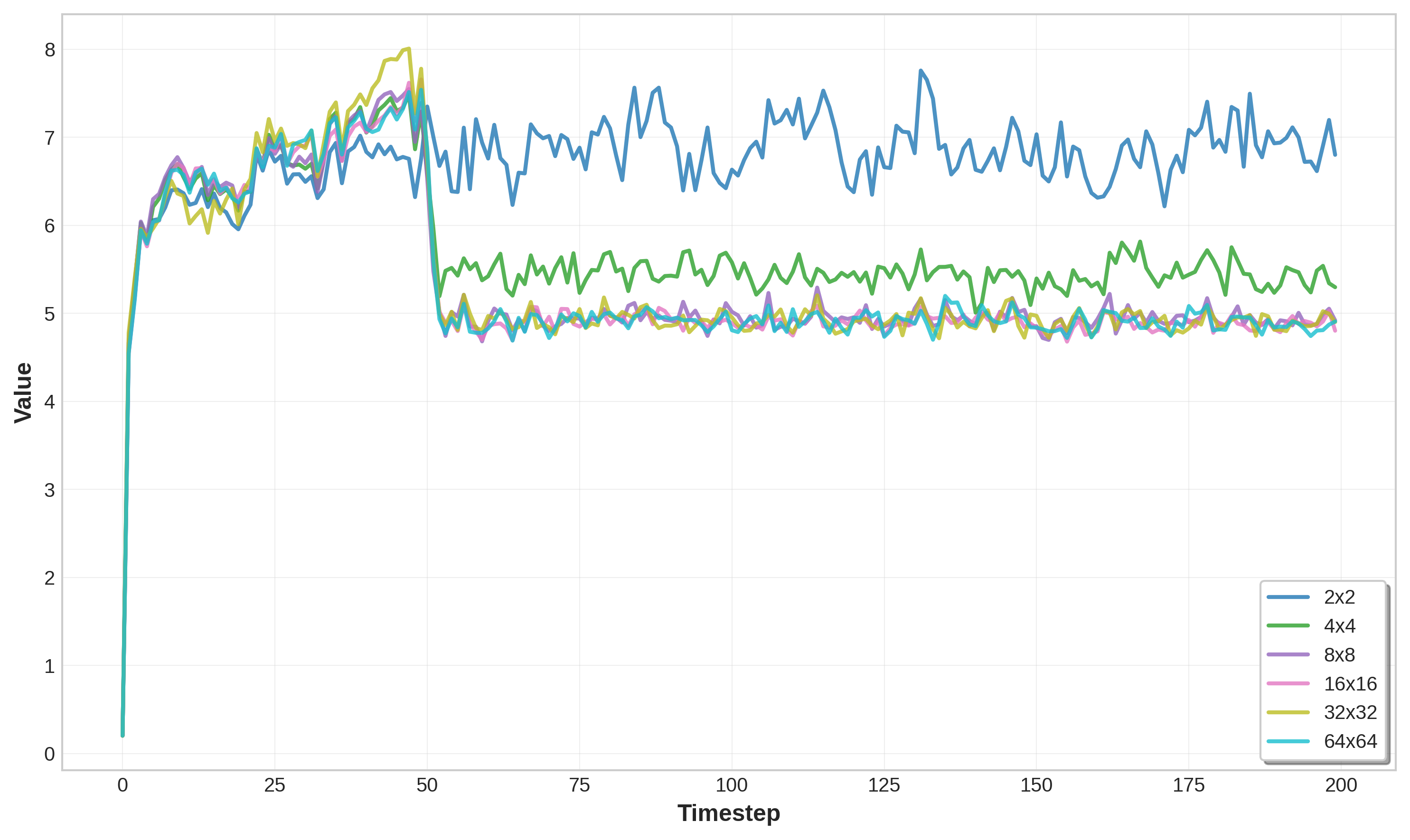}
         \caption{RMSE of particle‐filter estimates versus observation resolution (from $64\times64$ down to $2\times2$).}
         \label{fig: MSE_sum}
     \end{subfigure}
\end{figure}

Motivated by the application, we investigate the case where $\mcal{D} = \markspace \subset \R^2$. For computational reasons we choose to discretize the spatial domain into $4096$ sets, or, from an image analytical viewpoint, into $64\times 64$ pixels, whereas the decreasing observations' spatial resolutions are given as $64\times 64$, $32\times 32$, ..., $2 \times 2$ pixels; see Figure \ref{fig: resolution_figure_SPDE} for an example.

\paragraph{Particle filter estimations}Particle filters provide a numerical approximation of the Kushner-Stratonovich equation from Theorem \ref{kushner_strat_thm_multivariate}; see \cite[Ch. 8-10]{filtering_fundamentals} for details. A critical component of this approach involves calculating the forward steps of the Radon-Nikodym density $Z^M$, which, analogous to the signal and observation processes, is achieved using an explicit Euler scheme in our implemetation.

Let $ Y^M $ be a given observation according to \ref{introduction_MPP_discretized_observation_multivariate}. The \textit{ensemble size} \( L \in \N \) determines the the number of \textit{particles}, denoted by  \( X_{L,1}^M, \dots, X^M_{L,L} \), used in the particle filter. The algorithm iteratively simulates the particles' forward steps, assesses their likelihood, and then resamples them.
For a given time discretization $ t_1, \dots, t_N $ of $[0,T]$, the corresponding empirical distribution $ \frac{1}{L} \sum_{i=1}^L \delta_{X^M_{L,i}(t_j)} $ yields an approximation of the posterior distribution $ \eta^M_{t_j}$.

The empirical mean of the particles provides an estimate of the signal: 
\begin{equation} \label{simulations_particle_filter_est}
\overline{X}^M_L(t_j) := \frac{1}{L} \sum_{i=1}^L X^M_{L,i}(t_j) \approx X(t_j). 
\end{equation}
We assess the corresponding estimation error in \eqref{simulations_particle_filter_est} by computing the empirical mean squared errors.

The accuracy of \( \overline{X}^M_L \) depends on various factors, such as signal and observation noise amplitudes, the spatial resolution $M$, and the Monte Carlo sampling error, which decreases with larger ensemble sizes $L$. 
\paragraph{SPDE signal}
We consider a class of stochastic reaction-diffusion SPDEs, specifically of the form:
\begin{equation} \label{simulations_SPDE_dynamic}
    \begin{cases}
        \dif u(t) = \Big( \Delta u(t) + \varepsilon (u(t) - \alpha_1)(u(t) - \alpha_2)(\alpha_3 - u(t))  - v(t) + I \Big) \dif t+ B\dif W_1(t), \\
        \dif v(t) = \Big( \Delta v(t) + \gamma\big( \beta u(t) -v(t) \big) \Big) \dif t + \vartheta \dif W_2(t),
    \end{cases}
\end{equation}
which are commonly referred to as spatially extended stochastic FitzHugh-Nagumo dynamics, where $\Delta$ is the Neumann-Laplacian on $\mcal{D}$ and $W_1$ and $W_2$ are two independent cylindrical $Q$-Wiener processes.  In \cite{source_dicty_spde_1, source_dicty_spde_2} the stochastic FHN-System has been introduced as a spatially extendend stochastic two-phase dynamics to model and further analyze the actin dynamics in  D. discoideum. We set $X:=u$ as the signal process in our filtering problem, hence having an additional hidden process $v$ in the simulations. The parameters required to reproduce this experiment are available in the associated Github repository.

We applied a particle filter with $20.000$ particles to observations at various resolutions: \(64 \times 64\), \(32 \times 32\), down to \(2 \times 2\).
Using a spatial MSE computed for each time step , we observe that the estimation accuracy remains high even for relatively low-dimensional observations, as shown in Figure \ref{fig: MSE_sum}. One possible explanation is that, since the Laplacian is a is a strongly dissipative operator, its influence can still be captured effectively at lower resolutions, leading to accurate predictions of the signal state.

\subsection{Outlook: Filtering CLSM data}

The application of a Poisson particle filter to real CLSM data of \textit{D. discoideum} will be explored in future work. We plan to investigate parameter estimation under Poisson observation noise, expanding upon the theoretical framework established in \cite{toolbox_Gregor, source_dicty_spde_1}. While a detailed analysis will be provided in a forthcoming paper, we offer a brief overview of the intended applications.

Confocal Laser Scanning Microscopy (CLSM) is an optical imaging technique that enhances image clarity by selectively excluding out-of-focus light, effectively sectioning a three-dimensional object into thin two-dimensional optical slices. In CLSM, a laser beam is focused on single points within the sample, exciting fluorescent molecules that are present in the illuminated region. The sample is scanned point-by-point, and the emitted fluorescence passes through a pinhole aperture that blocks out-of-focus light, allowing only the fluorescence from the focal plane to reach the detector. This process results in an integer-valued photon count, which is typically transformed into a pixel value in a nonlinear fashion. In the analyzed data, we had access to the raw photon counts before their transformation into pixel values, allowing for more direct analysis of the imaging data.

\subsubsection{Data acquisition}

Experimental CLSM data was acquired using a laser scanning microscope (LSM780, Zeiss, Jena) equipped with a 20x objective lens and a 488 nm Argon laser. In order to access the raw photon counts, all recordings were performed under the "Photon Counting" acquisition mode. 

For the control experiments with fluorescein, a solution of 100 nM fluorescein sodium salt in Sørensen's buffer (14.7 mM KH$_2$PO$_4$, 2 mM Na$_2$HPO$_4$, pH 6.0)  was freshly prepared and further diluted to the desired final concentration before imaging. All fluorescein solutions were protected from light until imaging was performed. Timelapse recordings were acquired for 16 x 16 pixel frames, using a pixel dwell time of 16 $\mu$s, 40 $\mu$s or 81 $\mu$s, without any time delay between frames. 

For live cell imaging, we worked with giant \textit{D. discoideum} cells, produced through the electric pulse-induced fusion of individual cells \cite{source_giant_cells}. The cells (strain DdB NF1 KO, transformed with a plasmid for fluorescent labeling, SF108 as described in \cite{source_dicty_spde_2}) contain a green fluorescent protein that labels the intracellular actin (LifeAct-GFP). 
In all cases, samples were contained in a small petri dish with a glass bottom. 

\subsubsection{Poisson statistics in CLSM microscopy}
\begin{figure}
     \centering
     \caption{Photon count statistcs in CLSM data}
     \begin{subfigure}[t]{0.4\textwidth}
     \vskip 0pt
         \centering
         \includegraphics[width=\textwidth]{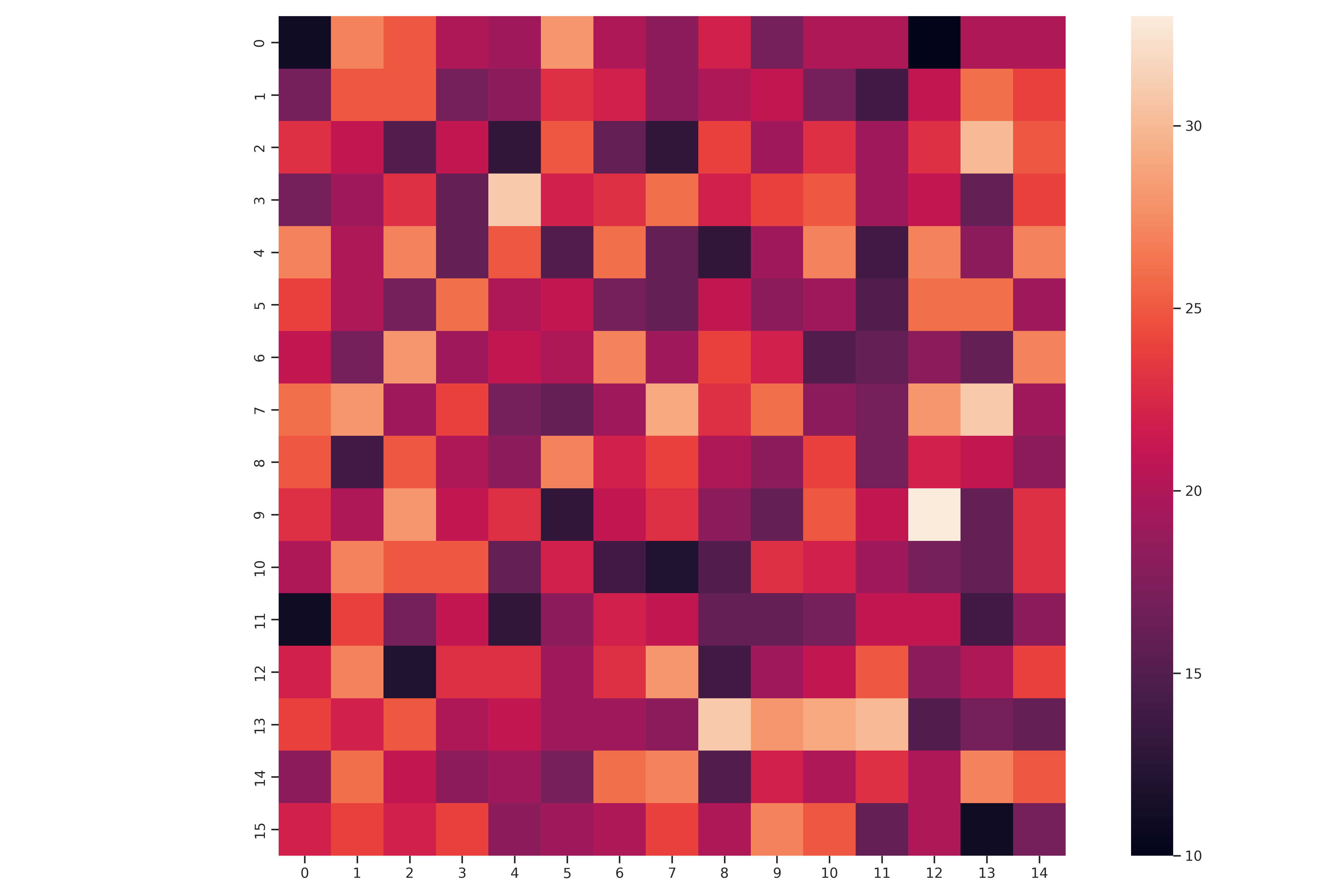} 
         \caption{A single $16 \times 16$ CLSM image from a representative dataset of 10,000 images, capturing a $10$ nM Fluorescein solution. }
         \label{fig:Fluorescin}
     \end{subfigure}
     \hfill
     \begin{subfigure}[t]{0.55\textwidth}
     \vskip 0pt
         \centering
         \includegraphics[width=\textwidth]{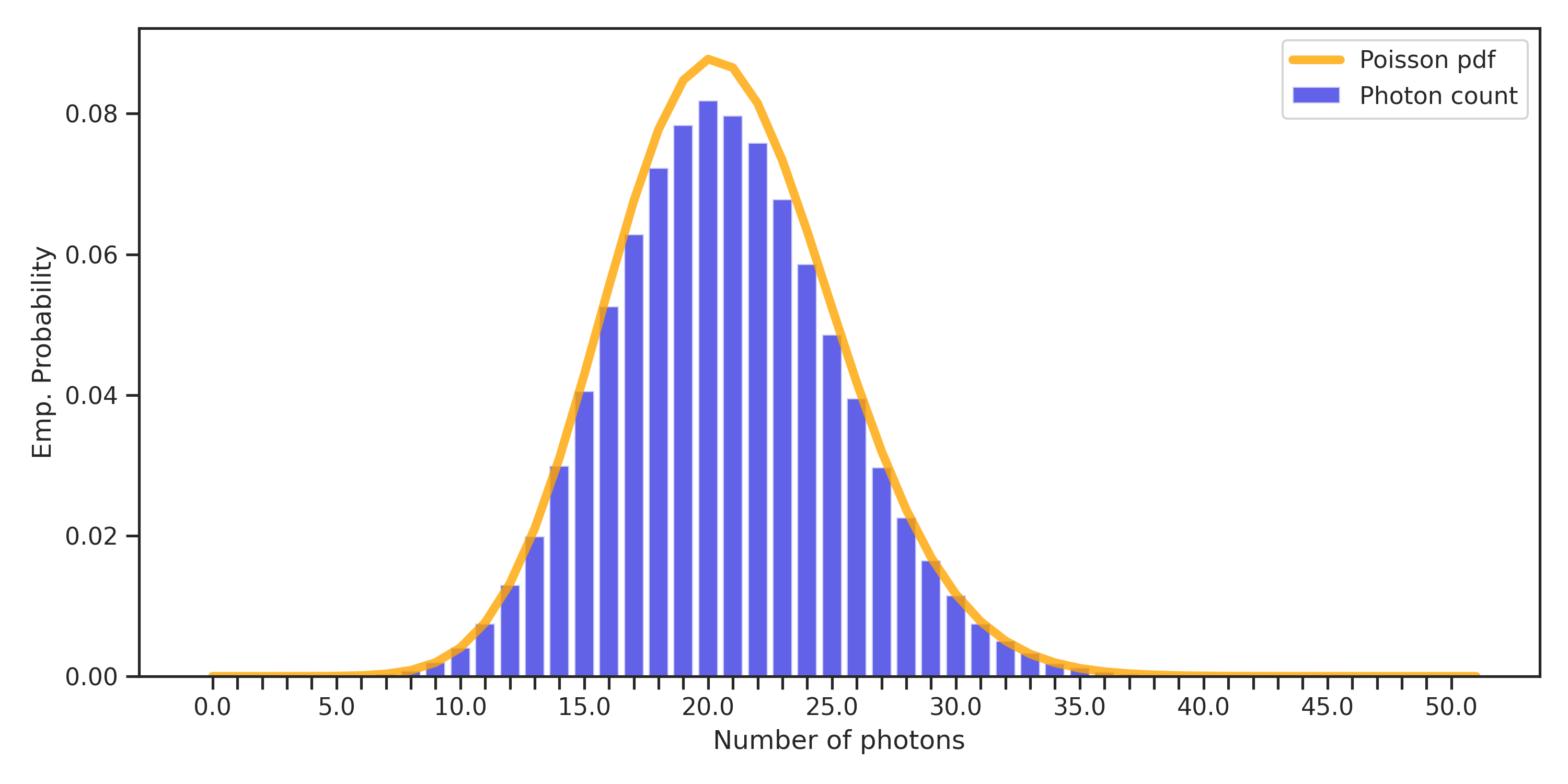} 
         \caption{Distribution of photon counts across all pixels from the entire set of images in the sample.}
         \label{fig:Photon statistics}
     \end{subfigure}
        \label{fig:Poisson noise CLSM}
\end{figure}

To validate the assumption that the observation noise in our data follows a Poisson distribution, commonly referred to as "shot noise" in statistical literature \cite{shot_noise_source}, we conducted an analysis on images of solutions containing varying concentrations of the fluorescent dye Fluorescein. Due to minimal diffusion over short time periods and within localized regions, it is reasonable to assume that the Fluorescein concentration remains approximately constant during the observations. An example of an image from such a "static" sample is shown in Figure \ref{fig:Fluorescin}. Each pixel in these images can be treated as a photon count sample from the same underlying Fluorescein concentration. We then compared the distribution of photon counts across all pixels with a Poisson probability density function (pdf) where the intensity parameter is given by the mean photon count, as illustrated in Figure \ref{fig:Photon statistics}. This analysis was performed across over 30 datasets, consistently showing that the bar plots of photon counts closely match Poisson distributions. The intensity of these distributions varied according to microscope settings, such as dwell time, laser intensity, and Fluorescein concentration. Further analysis revealed no significant correlation between photon counts, further supporting the Poisson noise assumption.

\subsubsection{Filtering CLSM data}
In a final experiment, we applied our filtering method to data obtained from confocal laser scanning microscopy recordings of giant \textit{D. discoideum} cells. Given that the datasets typically capture the entire cell, we began by extracting an area of interest (AOI) focused exclusively on the cell's interior to omit boundary effects \cite{source_dicty_spde_1}. The SPDE model \eqref{simulations_SPDE_dynamic} was used as the signal model, with parameters calibrated to ensure that $u$  maintains concentration values between $0$ and $1$ with large probability in good accordance with the observed data of actin concentrations. We assumed Poisson-distributed observation noise with an intensity of the form \eqref{simulations_intensity_form}, adjusting the scaling factor $c$ to align the model's photon counts with those observed in the data. 
\begin{figure}
    \caption{Poisson particle filter applied on real data}
     \centering
     \begin{subfigure}[b]{0.49\textwidth}
         \centering
         \includegraphics[width=\textwidth]{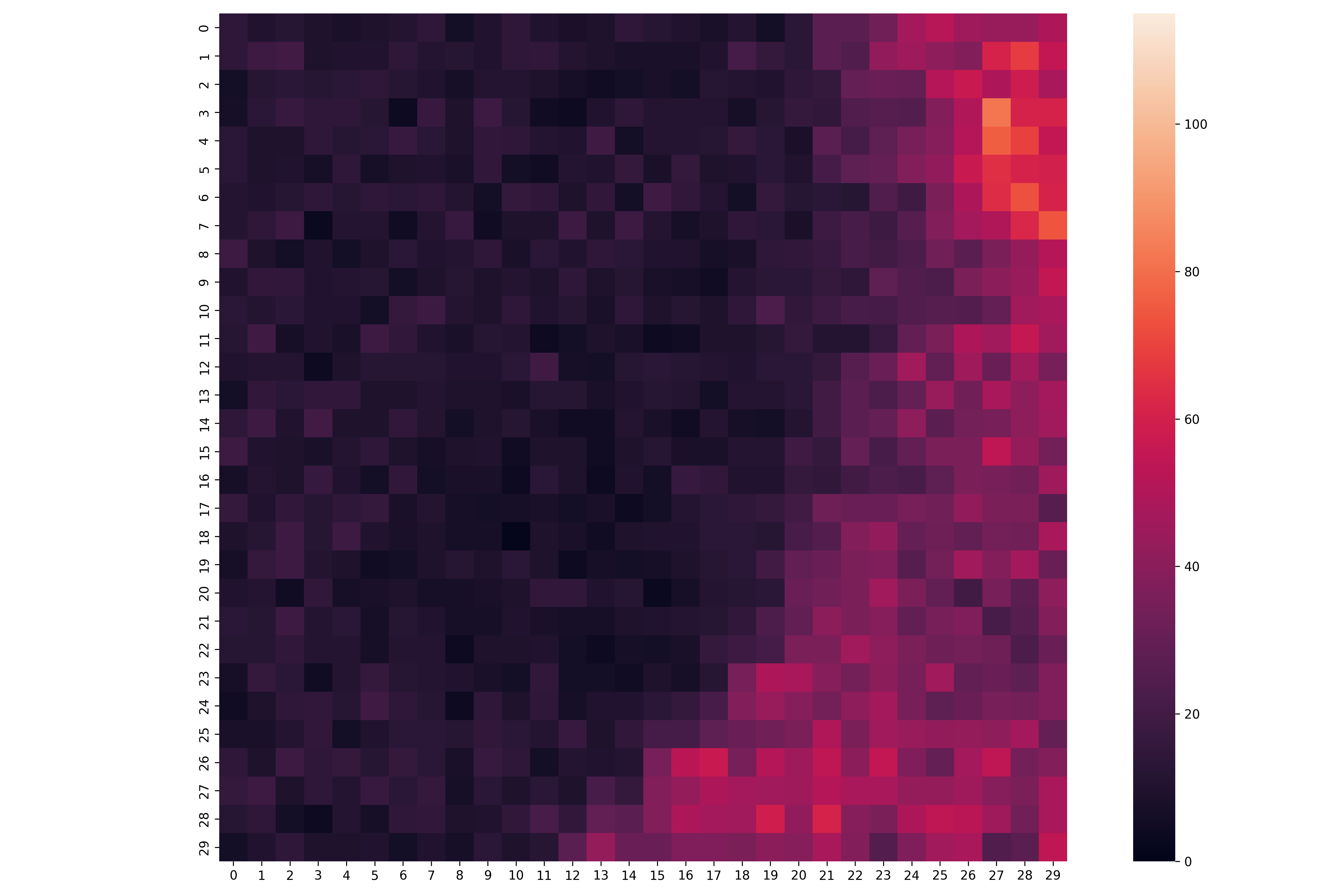}
         \caption{Extracted AOI of an CLSM image showing a giant \textit{D. discoideum cell}}
         \label{fig:CLSM real data}
     \end{subfigure}
     \hfill
     \begin{subfigure}[b]{0.49\textwidth}
         \centering
         \includegraphics[width=\textwidth]{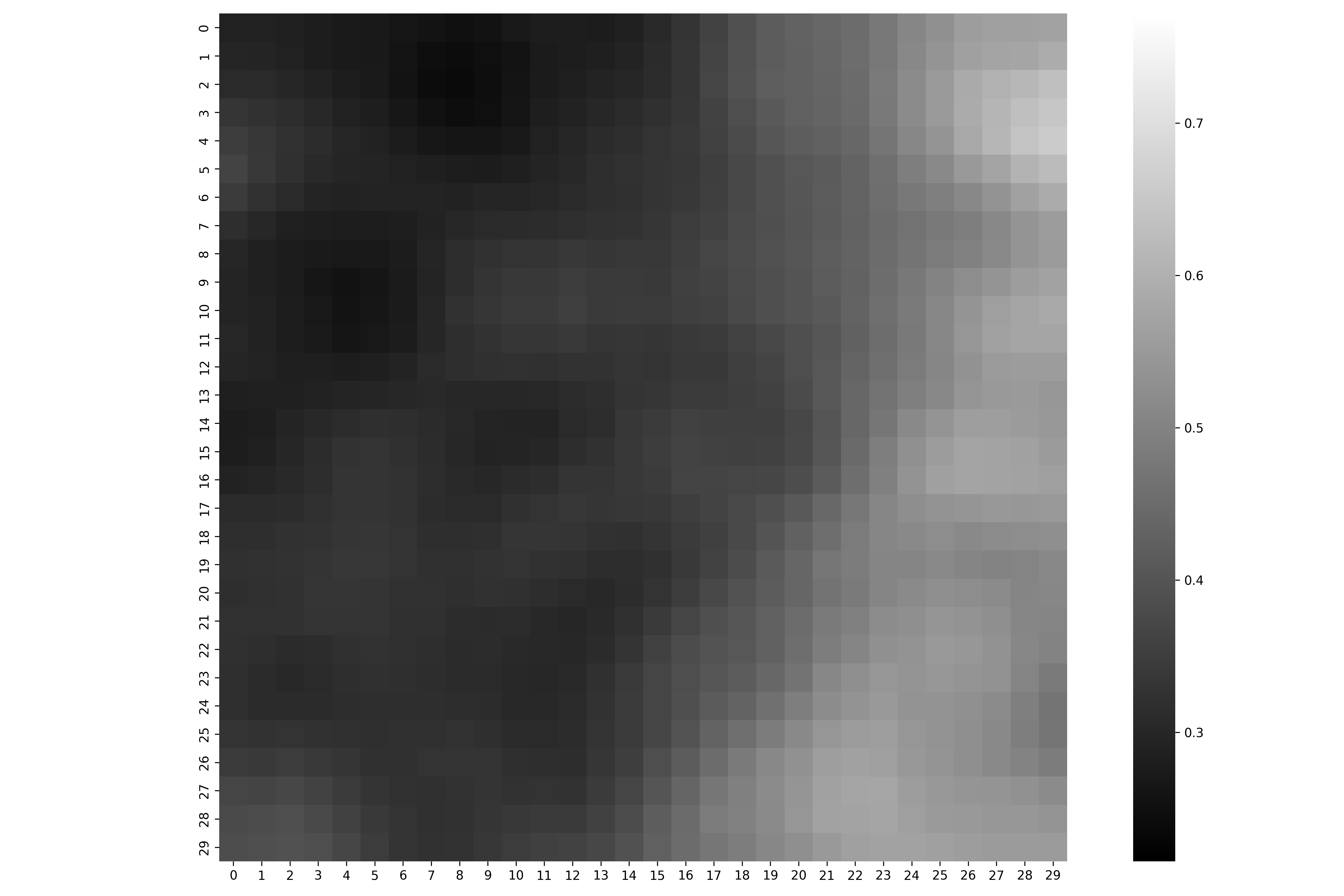} 
         \caption{Poisson particle filter estimation of the left image}
         \label{fig: CLSM estimaton}
     \end{subfigure}
        \label{fig:Filtering real data}
\end{figure}

Figure \ref{fig:Filtering real data} shows a data sample alongside the estimated state of the underlying actin dynamics. The experiments demonstrate that the filter effectively tracks wave-like actin movements, providing a satisfactory proof of concept across four different cell recordings.

While the initial results are promising, a significant challenge persists: the parameters must be manually selected, with no definitive method to ensure their accuracy beyond phenomenological validation. In future research, we aim to expand our theory and address this limitation by exploring parameter estimation techniques for SPDEs under point process noise, with a focus on potential applications in biophysics.

\subsection*{Acknowledgments}
This research has been partially funded by the  Deutsche Forschungsgemeinschaft (DFG)- Project-ID 318763901 - SFB1294, Project A01 “Statistics for Stochastic Partial Diﬀerential Equations” (JS and WS) and Project B02 "Inferring the dynamics underlying protrusion-driven cell motility" (CM-T).

\bibliography{main}
\appendix
\section*{Appendix}
\renewcommand{\thesection}{A} 
\setcounter{section}{1} 
\subsection*{Fundamentals of marked point processes}

 \paragraph{Simple point processes and MPPs} A point process $\chi$ on some state space $\mcal{S}$ is defined as  a measurable mapping from $(\Omega, \mcal{F}, \Prob)$ into $(\mcal{N}^\#_{\mcal{S}}, \Borel(\mcal{N}^\#_{\mcal{S}}))$, where $\mcal{N}^\#_\mcal{S}$ denotes the space of boundedly finite counting measures; see \cite[Ch. 9]{DVJ_vol2}  Motivated by our application, we choose $\mcal{S}:= [0,T] \times \markspace$ for $T$ from Section \ref{subsection_spde_solutions} and a compact set $\markspace \subset \R^{d_O}$, $d_O \in \N$. Let  $\mu_\markspace$ denote the $d_O$-dimensional Lebesgue measure. We introduce the measure space $(\markspace, \Borel(\markspace), \mu_\markspace)$, and call it the \textit{mark space} and are going to refer to $([0,T] \times \markspace, \Borel([0,T] \times \markspace), \dif t \times \mu_\markspace)$ when we speak of the \textit{product measure space}.
The following definitions and notations are taken from \cite[Ch.9]{DVJ_vol2}. 
\begin{definition} 
    \begin{enumerate}[(i)]
        \item By $\mcal{N}^{\#\ast}_{[0,T]}$ we denote the family of all \textit{simple} counting measures on $[0,T]$, meaning that for any $\zeta \in \mcal{N}^{\#\ast}_{[0,T]}$ we have
        \begin{equation}
            \zeta(\{t\}) \in \{0,1\} \text{ for all } t \in [0,T].
        \end{equation}
        \item By $\mcal{N}^{\#g}_{[0,T]\times \markspace}$ we denote the family of boundedly finite counting measures on the product measure space such that for any $\chi \in \mcal{N}^{\#g}_{[0,T]\times \markspace}$ the associated \textit{ground measure} $\chi^g$ defined by 
        \begin{equation}
            \chi^g(L) := \chi(L \times \markspace), \text{ for any } L \in \Borel([0,T]),
        \end{equation}
        is an element of $\mcal{N}^{\#\ast}_{[0,T]}$.
    \end{enumerate}

    \flushright $\diamond$
\end{definition}
Note that $\mcal{N}^{\#\ast}_{[0,T]}$ is not a closed subset of $\mcal{N}^{\#}_{[0,T]}$, and similarly, $\mcal{N}^{\#g}_{[0,T]\times \markspace}$ is not a closed subset of $\mcal{N}^{\#}_{[0,T]\times \markspace}$, as in general the existence of so-called accumulation points can not be ruled out.
Let $(\Omega, \mcal{F}, \Prob)$ be the filtered probability space from Section \ref{subsection_spde_solutions}.

\begin{definition} \label{definition_point_process}
    \begin{enumerate}[(i)]
        \item A \textit{point process} $\nu$ on the state space $[0,T]\times \markspace$ is a measurable mapping from $(\Omega, \mcal{F}, \Prob)$ into $(\mcal{N}^\#_{[0,T]\times \markspace}, \Borel(\mcal{N}^\#_{[0,T]\times \markspace}))$.
        \item A point process $\Bar{\nu} $ on $[0,T] $ is called \textit{simple} when $\Bar{\nu} \in  \mcal{N}^{\#\ast}_{[0,T]} \; \Prob$-a.s.
        \item A point process $\nu$ on $[0,T]\times \markspace$ is called \textit{marked point process (MPP)}  on $[0,T]$ with mark space $\markspace$ if $\nu \in \mcal{N}^{\#g}_{[0,T]\times \markspace} \; \Prob$-a.s.
        \item An MPP $\nu$ on $[0,T]\times \markspace$ is called \textit{marked Poisson process}  on $[0,T]$ with mark space $\markspace$ if its ground process is a Poisson process on $[0,T]$.
    \end{enumerate}
    \flushright $\diamond$
\end{definition}
As throughout the paper  $\markspace$ will always be the mark space, we are simply going to refer to any MPP  on $[0,T]$ with mark space $\markspace$ as an MPP on $[0,T]\times \markspace$. Sometimes it is also demanded that a marked Poisson process has a mark distribution which, given $\alpha$ is independent of $\lambda^g$; see \cite[p. 243]{Bremaud}. 

\begin{remark}[Finite boundedness on compact spaces]
    For any complete separable metric space $\mcal{S}$, denote by $\mcal{M}^\#_\mcal{S}$ the space of all boundedly finite measures on $\mcal{S}$, i.e., all countably additive, real-valued set functions $\xi$ with the property 
\begin{equation}
    \xi(A) < \infty \text{ for any bounded } A\in \Borel(\mcal{S}),
\end{equation}
and by  $\mcal{M}_\mcal{S}$ the family of all totally bounded measures on $\mcal{S}$. It is known that under the weak topology,  $\mcal{M}_\mcal{S}$ is complete separable metric space itself and that the family of all totally bounded counting measures $\mcal{N}_\mcal{S}$ is a closed subset of  $\mcal{M}_\mcal{S}$. Analogously, $\mcal{M}^\#_\mcal{S}$ is a complete separable metric space under the weak hash-topology, and the space of boundedly finite counting measures $\mcal{N}^\#_\mcal{S}$ is a closed subset of $\mcal{M}^\#_\mcal{S}$; see \cite[Ch. 9]{DVJ_vol2} for details. \newline

 It is evident that by compactness of $[0,T]\times\markspace$ the families  $\mcal{M}_{[0,T] \times \markspace}$ and $\mcal{M}^\#_{[0,T] \times \markspace}$, and thus also $\mcal{N}_{[0,T] \times \markspace}$ and $\mcal{N}^\#_{[0,T] \times \markspace}$, coincide. This implication will play a role in Section \ref{section_observation_schemes}, where we are going to exploit the fact that weak convergence on $\mcal{N}_{[0,T] \times \markspace}$ is metrizable to derive convergence rates; see \cite{Dudley_Baire}. However, keeping this identity in mind  we will stick to the notation using the $\#$-symbol for the measure spaces to be in line with point process literature.

\end{remark}

\paragraph{Doob-Meyer decomposition of MPPs}
For an MPP $\nu$, let us denote $\nu_\Gamma(t) := \nu((0,t] \times \Gamma)$ (and $\nu_\Gamma(0) :=  \nu(\{0\} \times \Gamma)$) for any $t \in [0,T]$ and $\Gamma \in \Borel(\markspace)$. Under mild assumptions, in particular boundedly finite first moment measure and absolute continuity of the so-called Campbell measure associated to $\nu$, see \cite[Ch. 13-14]{DVJ_vol2},  we have the existence of a $\Prob$-a.s. unique nonnegative \textit{conditional intensity} $\lambda$ w.r.t. $(\Prob, \mcal{F}_t)$, such that we have the integral representation
\begin{align} 
    \dif \nu_\Gamma(t) &=  \int_\Gamma \lambda(t, x) \; \mu_\markspace(\dif x) \, \dif t  + \dif N_\Gamma(t),
\end{align}
where the process $\tgeq{N_\Gamma(t)}$ defined by
\begin{equation} \label{sec_2_jump_martingale_general}
   \dif  N_\Gamma(t) := \dif \nu_\Gamma(t) -  \int_\Gamma \lambda(t, x) \; \mu_\markspace(\dif x) \, \dif t, \quad t \in (0,T],
\end{equation}
is a local right-continuous $\filtt$-martingale for any $\Gamma \in \Borel(\markspace)$. \newline

The analogous decomposition can be done for the ground measure $\nu^g$ of an MPP. There we simply introduce the \textit{ground process}  $\tgeq{\nu^g(t)}$ by 
\begin{equation}
    \nu^g(t) := \nu^g((0,t]) = \sum\limits_{(\tau_i, \kappa_i) \in \nu((0,t]\times \markspace)} \1\{(\tau_i, \kappa_i)  \in  (0,t] \times \markspace \}, \quad  t\in (0,T],
\end{equation}
which defines a right-continuous $\filtt$-adapted stochastic process. This leads to the integral representation
\begin{align} \label{semimartingale_decomposition_MPP_compensator}
    \dif \nu^g(t) &=  \int_\markspace \lambda(t, x) \; \mu_\markspace(\dif x) \, \dif t  + \dif N_\markspace(t).
\end{align}

It is often useful to factorize $\lambda$ into the intensity $\lambda^g$ of the ground process $Y^g$, defined $\Prob$-a.s. by
\begin{equation*}
    \lambda^g(t) := \int_\markspace \lambda(t,x) \, \mu_\markspace(\dif x), \quad t \in [0,T],  
\end{equation*}
and the stochastic kernel of the so-called \textit{conditional mark distribution} $\Phi(\dif x \cond t) :=  \phi(x \cond t) \mu_\markspace(\dif x)$ on $\markspace $, leading to the pair $\{\lambda^g(\cdot)\, , \, \Phi(\dif x \cond \cdot  ) \}$, called $(\Prob, \mcal{F}_t)$\textit{-local characteristics} in \cite{Bremaud}. The existence and uniqueness of such a factorization directly follows from the assumptions we made on the point process, see \cite[Prop. 14.3.II]{DVJ_vol2}. As they are derived directly from the compensator the conditional intensity, and equivalently the local characteristics, suffice to completely characterize an MPP w.r.t. $(\filtt)$.

\paragraph{Cox processes}
Finally,  all of the concepts in this paragraph can be easily extended to (marked) point processes, whose intensities $\lambda$ are  functions of some underlying random element $\zeta$. We provide the following definition and again refer to the standard books \cite{Bremaud, DVJ_vol1} for further details.
\begin{definition} Let $\zeta$ be a random measure on some measurable space $(S, \Borel(S))$.
    \begin{enumerate}[(i)]
        \item  An MPP $\nu$ on $[0,T]\times \markspace$ is a \textit{generalized marked Cox process} directed by $\zeta$, when its conditional intensity $\lambda$ is a measurable function of $\zeta$.
        \item  An MPP $\nu$ on $[0,T]\times \markspace$ is a \textit{marked Cox process} directed by $\zeta$, when it is a generalized marked Cox process whose ground process given $\zeta$ is a Poisson process on $[0,T]$; equivalently, given $\zeta$ the MPP $\nu$ is a marked Poisson process.
    \end{enumerate}
    \flushright $\diamond$
\end{definition}
We want to note that the notion of a generalized Cox process is not used consistently in the literature. In filtering theory it is standard procedure to let the random measure $\zeta$ be given as a nonnegative function  of the state $\xi(t)$ of some Markov process $\tgeq{\xi(t)}$. Equivalently, one can then say that the generalized Cox process is directed by $\tgeq{\xi(t)}$. An explicit construction will be given in the next section.

\end{document}